


\magnification 1200
\hsize 13.2cm
\vsize 20cm
\parskip 3pt plus 1pt
\parindent 5mm

\def\\{\hfil\break}


\font\seventeenbf=cmbx10 at 17.28pt

\font\twelvebf=cmbx10 at 12pt
\font\eightbf=cmbx8
\font\sixbf=cmbx6

\font\eighti=cmmi8
\font\sixi=cmmi6

\font\eightrm=cmr8
\font\sixrm=cmr6

\font\eightsy=cmsy8
\font\sixsy=cmsy6

\font\eightit=cmti8
\font\eighttt=cmtt8
\font\eightsl=cmsl8

\font\seventeenbsy=cmbsy10 at 17.28pt

\font\twelvebsy=cmbsy10 at 12pt
\font\tenbsy=cmbsy10
\font\eightbsy=cmbsy8
\font\sevenbsy=cmbsy7
\font\sixbsy=cmbsy6
\font\fivebsy=cmbsy5

\font\tenmsa=msam10

\font\sevenmsa=msam7
\font\fivemsa=msam5
\newfam\msafam
  \textfont\msafam=\tenmsa
  \scriptfont\msafam=\sevenmsa
  \scriptscriptfont\msafam=\fivemsa

\font\tenmsb=msbm10
\font\eightmsb=msbm8
\font\sevenmsb=msbm7
\font\fivemsb=msbm5
\newfam\msbfam
  \textfont\msbfam=\tenmsb
  \scriptfont\msbfam=\sevenmsb
  \scriptscriptfont\msbfam=\fivemsb
\def\Bbb{\fam\msbfam\tenmsb}

\font\tenCal=eusm10
\font\sevenCal=eusm7
\font\fiveCal=eusm5
\newfam\Calfam
  \textfont\Calfam=\tenCal
  \scriptfont\Calfam=\sevenCal
  \scriptscriptfont\Calfam=\fiveCal
\def\Cal{\fam\Calfam\tenCal}

\font\teneuf=eusm10
\font\teneuf=eufm10
\font\seveneuf=eufm7
\font\fiveeuf=eufm5
\newfam\euffam
  \textfont\euffam=\teneuf
  \scriptfont\euffam=\seveneuf
  \scriptscriptfont\euffam=\fiveeuf

\font\seventeenbfit=cmmib10 at 17.28pt

\font\twelvebfit=cmmib10 at 12pt
\font\tenbfit=cmmib10
\font\eightbfit=cmmib8
\font\sevenbfit=cmmib7
\font\sixbfit=cmmib6
\font\fivebfit=cmmib5
\newfam\bfitfam
  \textfont\bfitfam=\tenbfit
  \scriptfont\bfitfam=\sevenbfit
  \scriptscriptfont\bfitfam=\fivebfit


\catcode`\@=11
\def\eightpoint{%
  \textfont0=\eightrm \scriptfont0=\sixrm \scriptscriptfont0=\fiverm
  \def\rm{\fam\z@\eightrm}%
  \textfont1=\eighti \scriptfont1=\sixi \scriptscriptfont1=\fivei
  \def\oldstyle{\fam\@ne\eighti}%
  \textfont2=\eightsy \scriptfont2=\sixsy \scriptscriptfont2=\fivesy
  \textfont\itfam=\eightit
  \def\it{\fam\itfam\eightit}%
  \textfont\slfam=\eightsl
  \def\sl{\fam\slfam\eightsl}%
  \textfont\bffam=\eightbf \scriptfont\bffam=\sixbf
  \scriptscriptfont\bffam=\fivebf
  \def\bf{\fam\bffam\eightbf}%
  \textfont\ttfam=\eighttt
  \def\tt{\fam\ttfam\eighttt}%
  \textfont\msbfam=\eightmsb
  \def\Bbb{\fam\msbfam\eightmsb}%
  \abovedisplayskip=9pt plus 2pt minus 6pt
  \abovedisplayshortskip=0pt plus 2pt
  \belowdisplayskip=9pt plus 2pt minus 6pt
  \belowdisplayshortskip=5pt plus 2pt minus 3pt
  \smallskipamount=2pt plus 1pt minus 1pt
  \medskipamount=4pt plus 2pt minus 1pt
  \bigskipamount=9pt plus 3pt minus 3pt
  \normalbaselineskip=9pt
  \setbox\strutbox=\hbox{\vrule height7pt depth2pt width0pt}%
  \let\bigf@ntpc=\eightrm \let\smallf@ntpc=\sixrm
  \normalbaselines\rm}
\catcode`\@=12

\def\eightpointbf{%
 \textfont0=\eightbf   \scriptfont0=\sixbf   \scriptscriptfont0=\fivebf
 \textfont1=\eightbfit \scriptfont1=\sixbfit \scriptscriptfont1=\fivebfit
 \textfont2=\eightbsy  \scriptfont2=\sixbsy  \scriptscriptfont2=\fivebsy
 \eightbf
 \baselineskip=10pt}

\def\tenpointbf{%
 \textfont0=\tenbf   \scriptfont0=\sevenbf   \scriptscriptfont0=\fivebf
 \textfont1=\tenbfit \scriptfont1=\sevenbfit \scriptscriptfont1=\fivebfit
 \textfont2=\tenbsy  \scriptfont2=\sevenbsy  \scriptscriptfont2=\fivebsy
 \tenbf}

\def\twelvepointbf{%
 \textfont0=\twelvebf   \scriptfont0=\eightbf   \scriptscriptfont0=\sixbf
 \textfont1=\twelvebfit \scriptfont1=\eightbfit \scriptscriptfont1=\sixbfit
 \textfont2=\twelvebsy  \scriptfont2=\eightbsy  \scriptscriptfont2=\sixbsy
 \twelvebf
 \baselineskip=14.4pt}

\def\seventeenpointbf{%
 \textfont0=\seventeenbf  \scriptfont0=\twelvebf  \scriptscriptfont0=\eightbf
 \textfont1=\seventeenbfit\scriptfont1=\twelvebfit\scriptscriptfont1=\eightbfit
 \textfont2=\seventeenbsy \scriptfont2=\twelvebsy \scriptscriptfont2=\eightbsy
 \seventeenbf
 \baselineskip=20.736pt}


\newdimen\srdim \srdim=\hsize
\newdimen\irdim \irdim=\hsize
\def\NOSECTREF#1{\noindent\hbox to \srdim{\null\dotfill ???(#1)}}
\def\SECTREF#1{\noindent\hbox to \srdim{\csname REF\romannumeral#1\endcsname}}
\def\INDREF#1{\noindent\hbox to \irdim{\csname IND\romannumeral#1\endcsname}}
\newlinechar=`\^^J
\def\openauxfile{
  \immediate\openin1\jobname.aux
  \ifeof1
  \message{^^JCAUTION\string: you MUST run TeX a second time^^J}
  \let\sectref=\NOSECTREF \let\indref=\NOSECTREF
  \else
  \input \jobname.aux
  \message{^^JCAUTION\string: if the file has just been modified you may
    have to run TeX twice^^J}
  \let\sectref=\SECTREF \let\indref=\INDREF
  \fi
  \message{to get correct page numbers displayed in Contents or Index
    Tables^^J}
  \immediate\openout1=\jobname.aux
  \let\END=\end \def\end{\immediate\closeout1\END}}

\newbox\titlebox   \setbox\titlebox\hbox{\hfil}
\newbox\sectionbox \setbox\sectionbox\hbox{\hfil}
\def\folio{\ifnum\pageno=1 \hfil \else \ifodd\pageno
           \hfil {\eightpoint\copy\sectionbox\kern8mm\number\pageno}\else
           {\eightpoint\number\pageno\kern8mm\copy\titlebox}\hfil \fi\fi}
\footline={\hfil}
\headline={\folio}

\def\titlerunning#1{\setbox\titlebox\hbox{\eightpoint #1}}
\def\title#1{\noindent\hfil$\smash{\hbox{\seventeenpointbf #1}}$\hfil
             \titlerunning{#1}\medskip}

\newcount\numbersection \numbersection=-1
\def\sectionrunning#1{\setbox\sectionbox\hbox{\eightpoint #1}
  \immediate\write1{\string\def \string\REF
      \romannumeral\numbersection \string{%
      \noexpand#1 \string\dotfill \space \number\pageno \string}}}
\def\section#1{%
  \par\vskip0.666cm\penalty -100
  \vbox{\baselineskip=14.4pt\noindent{{\twelvepointbf #1}}}
  \vskip2pt
  \penalty 500
  \advance\numbersection by 1
  \sectionrunning{#1}}

\def\subsection#1{%
  \par\vskip0.5cm\penalty -100
  \vbox{\noindent{{\tenpointbf #1}}}
  \vskip1pt
  \penalty 500}

\newcount\numberindex \numberindex=0
\def\index#1#2{%
  \advance\numberindex by 1
  \immediate\write1{\string\def \string\IND #1%
     \romannumeral\numberindex \string{%
     \noexpand#2 \string\dotfill \space \string\S \number\numbersection,
     p.\string\ \space\number\pageno \string}}}

\newdimen\itemindent \itemindent=\parindent

\def\item#1{\par\noindent\hangindent\itemindent%
            \rlap{#1}\kern\itemindent\ignorespaces}
\def\itemitem#1{\par\noindent\hangindent2\itemindent%
            \kern\itemindent\rlap{#1}\kern\itemindent\ignorespaces}
\def\itemitemitem#1{\par\noindent\hangindent3\itemindent%
            \kern2\itemindent\rlap{#1}\kern\itemindent\ignorespaces}

\long\def\claim#1|#2\endclaim{\par\vskip 5pt\noindent
{\tenpointbf #1.}\ {\it #2}\par\vskip 5pt}

\def\proof{\noindent{\it Proof}}

\def\today{\ifcase\month\or
January\or February\or March\or April\or May\or June\or July\or August\or
September\or October\or November\or December\fi \space\number\day,
\number\year}

\catcode`\@=11
\newcount\@tempcnta \newcount\@tempcntb
\def\timeofday{{%
\@tempcnta=\time \divide\@tempcnta by 60 \@tempcntb=\@tempcnta
\multiply\@tempcntb by -60 \advance\@tempcntb by \time
\ifnum\@tempcntb > 9 \number\@tempcnta:\number\@tempcntb
  \else\number\@tempcnta:0\number\@tempcntb\fi}}
\catcode`\@=12

\def\bibitem#1&#2&#3&#4&%
{\hangindent=0.8cm\hangafter=1
\noindent\rlap{\hbox{\eightpointbf #1}}\kern0.8cm{\rm #2}{\it #3}{\rm #4.}}


\def\bD{{\Bbb D}}

\def\bQ{{\Bbb Q}}
\def\bR{{\Bbb R}}

\def\bZ{{\Bbb Z}}


\def\cE{{\Cal E}}

\def\cI{{\Cal I}}

\def\cO{{\Cal O}}
\def\cR{{\Cal R}}

\def\cT{{\Cal T}}
\def\cX{{\Cal X}}
\def\cP{{\Cal P}}
\def\cV{{\Cal V}}


\def\bu{{\scriptstyle\bullet}}

\def\square{{\hfill \hbox{
\vrule height 1.453ex  width 0.093ex  depth 0ex
\vrule height 1.5ex  width 1.3ex  depth -1.407ex\kern-0.1ex
\vrule height 1.453ex  width 0.093ex  depth 0ex\kern-1.35ex
\vrule height 0.093ex  width 1.3ex  depth 0ex}}}
\def\qed{\kern10pt$\square$}
\def\hexnbr#1{\ifnum#1<10 \number#1\else
 \ifnum#1=10 A\else\ifnum#1=11 B\else\ifnum#1=12 C\else
 \ifnum#1=13 D\else\ifnum#1=14 E\else\ifnum#1=15 F\fi\fi\fi\fi\fi\fi\fi}
\def\msatype{\hexnbr\msafam}
\def\msbtype{\hexnbr\msbfam}
\mathchardef\restriction="3\msatype16   
\mathchardef\boxsquare="3\msatype03
\mathchardef\preccurlyeq="3\msatype34
\mathchardef\compact="3\msatype62
\mathchardef\smallsetminus="2\msbtype72   
\mathchardef\subsetneq="3\msbtype28
\mathchardef\supsetneq="3\msbtype29
\mathchardef\leqslant="3\msatype36   
\mathchardef\geqslant="3\msatype3E   
\mathchardef\stimes="2\msatype02
\mathchardef\ltimes="2\msbtype6E
\mathchardef\rtimes="2\msbtype6F

\def\ddbar{\partial\overline\partial}

\let\ol=\overline

\let\wt=\widetilde
\let\wh=\widehat
\let\text=\hbox
\def\buildo#1^#2{\mathop{#1}\limits^{#2}}
\def\buildu#1_#2{\mathop{#1}\limits_{#2}}
\def\ort{\mathop{\hbox{\kern1pt\vrule width4.0pt height0.4pt depth0pt
    \vrule width0.4pt height6.0pt depth0pt\kern3.5pt}}}

\def\vlra{\mathrel{\smash-}\joinrel\mathrel{\smash-}\joinrel%
    \kern-2pt\longrightarrow}
\def\srelbar{\vrule width0.6ex height0.65ex depth-0.55ex}
\def\merto{\mathrel{\srelbar\kern1.3pt\srelbar\kern1.3pt\srelbar
    \kern1.3pt\srelbar\kern-1ex\raise0.28ex\hbox{${\scriptscriptstyle>}$}}}




\long\def\InsertFig#1 #2 #3 #4\EndFig{\par
\hbox{\hskip #1mm$\vbox to#2mm{\vfil\special{"
(/home/demailly/psinputs/grlib.ps) run
#3}}#4$}}
\long\def\LabelTeX#1 #2 #3\ELTX{\rlap{\kern#1mm\raise#2mm\hbox{#3}}}


\itemindent = 7mm

\title{Quantitative extensions of pluricanonical }
\title{forms and closed positive currents}
\titlerunning{}

\vskip10pt

\centerline {{\tenrm Bo BERNDTSSON} and  {\tenrm Mihai P\u AUN}}

\vskip20pt

\noindent{\bf Abstract. \it {In this article we establish several 
Ohsawa-Takegoshi type theorems for twisted pluricanonical forms and metrics of adjoint $\bR$-bundles.}}

%
\section{\S0 Introduction}

\medskip 
\noindent Our main goal in this article is to generalize the Ohsawa-Takegoshi extension theorem in the context of pluricanonical forms. 

Let ${\cX}\to \Delta$ be a smooth projective family, and let
$(L, h_L)\to \cX$ be a line bundle endowed with a metric $h_L$.
\noindent The ''standard`` assumptions for the metric $h_L$ are:
{\itemindent 3mm
\item {$\bullet$}
The curvature current of $(L,h_L)$ is positive, i.e. $\sqrt {-1}\Theta_{h_L}(L)\geq 0$;
\smallskip
\item{$\bullet$}
 The restriction of the metric $h_L$ to the central fiber is well defined 
 $h_{L\vert {\cX}_0}\not\equiv\infty$.
}
\smallskip
\noindent Under these circumstances, the extension theorem established in  [32] (and subsequently developed in [1], [2], [13], [27], [33], [34], [41], [46]) states as follows:
{\sl let $u$ be a holomorphic section of the bundle $\displaystyle K_{\cX}+ L_{|\cX_0}$ which is $L^2$ with respect to $h_L$, i.e.  
$$\int_{\cX_0}|u|^2e^{-\varphi_L}< \infty.\leqno(1)$$
Then there exists a section $U$ of the bundle 
$\displaystyle K_{\cX}+ L$ whose restriction to $\cX_0$ is equal to $u$, and such that 
$$\int_{\cX}|U|^2e^{-\varphi_L}\leq C_0\int_{\cX_0}|u|^2e^{-\varphi_L}
$$
where $C_0$ is a purely numerical constant.} The meaning of the word 
``restriction" above is that over the central fiber we have $\displaystyle U_{|\cX_0}= u\wedge dt$.

\noindent 
We prove in this article similar effective extension statements for bundles of type
$$pK_{\cX}+ L.\leqno(2)$$
If $p\geq 2$, a first result to be mentioned is 
the {\sl invariance of plurigenera} due to Y.-T. Siu (cf. [41]), which completely elucidates the case of pluricanonical forms (i.e. without the additional twisting $L$). Indeed, it is possible (and not very complicated) to refine further the result in [41], and obtain the following statement:
{\sl let $u$ be a holomorphic section of the bundle 
$\displaystyle pK_{\cX_0}$;
then there exists a section $U$ of the bundle 
$\displaystyle pK_{\cX}$ whose restriction to $\cX_0$ is equal to $u$, and such that 
$$\int_{\cX}|U|^{2\over p}\leq C_0\int_{\cX_0}|u|^{2\over p}$$
where $C_0$ is the same constant as above.}

Motivated by applications in algebraic geometry, one has to generalize this kind of results for {\sl twisted pluricanonical forms}
but in this setting, the {\sl optimal} integrability conditions to be imposed are less clear: for example, replacing (1) with the natural $L^{2\over p}$ convergence is not enough (cf. the examples in [16], [21]). 
\medskip 

\noindent We describe next the results we obtain in the present article.

\smallskip

To start with, we recall the following notion. Consider an ideal ${\cI}\subset {\cO}_X$ and a positive integer $k\geq 0$. We denote by 
$\ol {\cI}^{(k)}$ the integral closure of the $k^{\rm th}$ power of $\cI$: it is the 
ideal constructed as follows. Let $x\in X$ be an arbitrary point, and let $(f_1,..., f_r)$
be the generators of ${\cI}$ at $(X, x)$. Then $\ol {\cI}^{(k)}$ is described locally at $x$ as follows
$$\ol {\cI}^{(k)}_x:= \big\{g\in {\cO}_{X, x}/ |g|^{2\over k}\leq \sum_\alpha |f_\alpha|^2\big\}.$$

\vskip 7pt 
\noindent In this context, we first prove the next refined version of the 
twisted invariance of plurigenera.

\claim 0.1 Theorem|Let $\pi: \cX\to \bD$ be a projective
family over the unit disk and
let $(L, h)$ be a
hermitian line bundle, with the properties $\bullet$ and $\bullet$ above. Then there exists a universal constant 
$C_0> 0$ such that for any positive integers $p\geq q$ and for any section $$u\in \displaystyle H^0\Bigl(\cX_0, pK_{\cX_0}+ qL)\otimes 
\ol{\cI}
(h_{L\vert \cX_0})^{(q)}\Bigr)$$ 
there exists a section $$U\in \displaystyle  
H^0\bigl(\cX, pK_{\cX}+ qL\bigr)$$ 
such that: 
{\itemindent 6mm
\item {(i)} Over the central fiber we have $U_{|\cX_0}= u\otimes d\pi^{\otimes p}$;
\smallskip
\item {(ii)} The next
$L^{2/p}$ integrability condition holds
$$\int_{\cX}|U|^{2\over p} e^{-{q\over p}\varphi_L}\leq C_0 \int_{\cX_0}
|u|^{2\over p} e^{-{q\over p}\varphi_L}.
$$}
\endclaim

\vskip 7pt

If $p=q= 1$, then this is precisely the Ohsawa-Takegoshi extension theorem recalled above. 
For $p\geq 2$ the origins of the qualitative part of our result is the work of Siu
see [40], [41], and also [10], [14], [16], [17], [18], [19], [20], [21], [25], 
[26], [27], [29], [30], [35], [44], [45], [46], [47], [49] for related statements.
\vskip 5pt
\noindent An important source of inspiration for theorem 0.1 arise in particular from the results obtained by 
H. Tsuji, S.~Takayama and
C. Hacon-J. McKernan respectively in connection with their work on pluricanonical series (see [46], [44], [20]). 
To make this more transparent, we consider the following variant of 0.1.
\claim 0.1$^\prime$ Theorem|Let $X$  be projective, and let $(L, h_L)$ be a hermitian line bundle on $X$. Let 
$S\subset X$ be a non-singular, irreducible submanifold of codimension 1, such that 
 $\displaystyle h_{\vert S}$ is well defined. Assume that the next curvature condition
 is satisfied
$$\Theta_{h_L}(L)\geq \varepsilon_0\omega.$$

\noindent Then any section of $\displaystyle (pK_S+ qL)\otimes \ol {\cI}
(h_{L\vert S})^{(q)}$ extends to $X$ as a section of the line bundle
$p(K_X+ S) +qL.$ 
\endclaim

\noindent This result is not a consequence of theorem 0.1, but of its proof: the only difference is
the version of the Ohsawa-Takegoshi theorem to be used in the inductive process.
Moreover, following [27] it is possible to formulate (and 
prove) the higher codimensional analogue of the previous statement
--where the hypersurface $S$ will be replaced by a {\sl maximal center} of some $\bQ$--divisor. \hfill\qed

\medskip


We remark that in theorem 0.1 it is not enough to assume the convergence of the integral
$\displaystyle \int_{\cX_0}
|u|^{2\over p} e^{-{q\over p}\varphi_L}$ in order to infer the extension of $u$; the additional hypothesis 
$$u\in \ol{\cI}
(h_{L\vert \cX_0})^{(q)}$$
is 
needed. For many purposes however it is desirable the avoid this latter condition, i.e. to replace it with something more ``manageable". This is the problem we address in the second part of our article,
where we assume for simplicity that $q= 1$ and $\displaystyle \cI(h_{L|\cX_0}^{{1/p}})= \cO_{\cX_0}$. We notice that in general this latter hypothesis does not implies that $\displaystyle \cI(h_{L|\cX_0}^{{}})= \cO_{\cX_0}$.

We will analyze here the extension of sections of 
(2) under the hypothesis that the curvature current of $L$ is only assumed to be semi-positive. Hence, unlike the usual setting, the bundle $L$ or its restriction to the central fiber is not necessarily big, but a {\sl natural vanishing assumption} for the section to be extended is needed. Our next result
can be seen as an effective version of the Ein-Popa theorem in [16]; also, it is a generalization of results due to J.-P. Demailly and H. Tsuji in [15], respectively [47], [48].
There are many notations/hypothesis we have to introduce before stating it, but they are natural in the context of the study of twisted 
pluricanonical systems on algebraic manifolds.

Let $\pi:\cX\to \bD$ be a proper, surjective map, where $\bD$ is the unit disk.
We assume that the central fiber $\cX_0= \pi^{-1}(0)$ is non-singular, and let $L\to \cX$
be a hermitian line bundle such that $c_1(L)$ contains the current
$$p([\Delta]+ \alpha)\in c_1(L)\leqno(3)$$
where the notations are as follows.

\noindent ${\bf (a)}$ $\Delta:= \sum _{j\in J}\nu^jY_j$ is an {\sl effective} $\bQ$-divisor, such that $p\nu^j\in \bZ$ for any $j\in J$; the hypersurfaces 
$Y_j\subset \cX$ together with $\cX_0$ have strictly normal crossings.

\noindent ${\bf (b)}$ $\alpha$ is a closed, non-singular, semi-positive form
of (1,1)-type, with the property that $\{p\alpha\}\in H^2(\cX, \bZ)$.
\smallskip

\noindent Furthermore, we assume that the bundle $K_{\cX}+ 1/pL$
is pseudoeffective, and let $h_{\rm min}$ be a metric with minimal singularities corresponding to it; we denote by $\Theta_{\rm min}$ its curvature current. We assume that 
$$\nu_{\rm min}\big(\{K_{\cX}+ 1/pL\}, \cX_0\big)= 0\leqno(4)$$
that is to say, the minimal multiplicity of the class $\{K_{\cX}+ 1/pL\}$ 
along the
central fiber $\cX_0$ is equal to zero (see e.g. [6]). 
Let $A\to \cX$ be an ample line bundle. The assumption (4) implies that the metric with minimal singularities $h_{\rm min, \varepsilon}$ corresponding to the class
$K_{\cX}+ 1/pL+ \varepsilon A$ is not identically $+\infty$ when restricted to 
$\cX_0$ (see [6]),  so that we can write
$$\Theta_{\rm min, \varepsilon|\cX_0}= \sum_{j\in J}\rho^j_{\rm min, \varepsilon}[Y_{j0}]+ \Lambda_{0\varepsilon}\leqno(5)$$
where $Y_{j0}:= Y_j\cap \cX_0$ and where $(\rho^j_{\rm min, \varepsilon})$ are 
positive real numbers. For each $j$, the sequence 
$(\rho^j_{\rm min, \varepsilon})$ is decreasing, and we define
$$\rho^j_{\rm min, \infty}:= \lim_{\varepsilon\to 0}\rho^j_{\rm min, \varepsilon}.\leqno(6)$$ 
We introduce the notation
$$J^\prime:= \{j\in J: \rho^j_{\rm min, \infty}< \nu^j\}$$
and we assume furthermore that the next condition is satisfied:

\noindent ${\bf (c)}$ We have $\nu^j\leq 1$ and
for any subset $I\subset J^\prime$ and any $\varepsilon> 0$ 
the restriction of the current 
$\Lambda_{0 \varepsilon}$ defined in (6) to the intersection $\displaystyle \bigcap_{m\in I}Y_m\cap\cX_0$
is well-defined.

\smallskip

\noindent Let $h_0= e^{-\varphi_0}$ be a metric on the $\bQ$-bundle
$K_{\cX_0}+ 1/pL$ with the property that
$$\Theta_{h_0}(K_{\cX_0}+ 1/pL)\geq 0$$
and such that the following inequalities are satisfied
$$\varphi_0\leq \sum_{j\in J^\prime}\rho^j_{\rm min, \infty}\log|f_{Y_j}|^2+ 
\sum_{j\in J\setminus J^\prime}\nu^j\log|f_{Y_j}|^2\leqno(7)$$
and
$$\int_{\cX_0}e^{\varphi_0-{1\over p}\varphi_L}< \infty\leqno(8)$$
We denote by $\varphi_L$ the singular metric on $L$ induced by the 
decomposition (3), and 
for each $j$ we denote by $f_{Y_j}$ in (7)
the local equations of the 
hypersurface $Y_j$. We state now our next result.
\claim 0.2 Theorem|Under the hypothesis $\bf(a)-(c)$ and $\rm (7), (8)$
above, the restriction $\varphi_{\rm min|\cX_0}$ is well-defined, and there exists a constant $C< 0$ depending only on the quantity $\rm (8)$ and the geometry of the map $\pi$ such that the following inequality holds at each point of $\cX_0$ 
$$\varphi_{\rm min|\cX_0}\geq C+ \varphi_0.\leqno(9)$$
If we assume moreover that $\nu^j< 1$ for all $j$, then given any section $u$ of the bundle $pK_{\cX_0}+ L$
whose zero divisor is greater than 
$$p\sum_{j\in J^\prime}\rho^j_{\rm min, \infty}[Y_{j0}]+ 
p\sum_{j\in J\setminus J^\prime}\nu^j[Y_{j0}]\leqno(\star)$$
there exists a section $U$ of $pK_{\cX}+ L$ extending $u$, and such that
$$\int_{\cX}|U|^{2\over p} e^{-{1\over p}\varphi_L}\leq C_0 \int_{\cX_0}
|u|^{2\over p} e^{-{1\over p}\varphi_L}.
$$
\hfill\qed
\endclaim
\noindent  Actually, we obtain an even 
more precise result:  we can replace the metric $\varphi_{\rm min}$ in (9)
with the so-called {\sl super-canonical metric} on the bundle 
$\displaystyle K_{\cX/\bD}+ {1\over p}L$ (in the terminology of [15], [47]). We prefer however the formulation above, for reasons that will appear in a moment.

We remark that as a consequence of (9) we obtain Ohsawa-Takegoshi type estimates for the extension $U$, provided that the section $u$ vanishes along the divisor 
$(\star)$.

If the form $\alpha$ in $\bf (b)$ is strictly positive, then the second part of the preceding result
was established in [22], [16]. Also, we refer to the section 17 of the article [15] (and the references therein) for an enlightening introduction and related results around this circle of ideas.
\hfill\qed
\medskip 

\noindent In order to give another interpretation of the result 0.2, we assume that we have $\nu^j< 1$, i.e. the pair $(\cX, \Delta)$ is {\sl klt} in algebro-geometric
language. 

Let $L^\prime= L_{|\cX_0}- p\sum_{j\in J^\prime}\rho^j_{\rm min, \infty}[Y_{j0}]-
p\sum_{j\in J\setminus J^\prime}\nu^j[Y_{j0}]$; it is not too difficult to show that the 
bundle $K_{\cX_0}+ 1/pL^\prime$ is pseudoeffective (see e.g. the arguments at the end of section $B$). 
We denote by $\varphi_{\rm min}^\prime$ the metric with minimal singularities corresponding to the bundle $K_{\cX_0}+ 1/pL^\prime$; then we have
$$\big|\varphi_{\rm min|\cX_0}- \sum_{j\in J^\prime}\rho^j_{\rm min, \infty}\log|f_j|^2- 
\sum_{j\in J\setminus J^\prime}\nu^j\log|f_j|^2 - \varphi_{\rm min}^\prime\big|\leq C.\leqno (10)$$
so the singularities of the restriction $\varphi_{\rm min|\cX_0}$ are completely understood in terms of the extremal metric $\varphi_{\rm min}^\prime$. Except for the {\sl rationality} of the coefficients $\rho^j_{\rm min, \infty}$, the relation (10) is the metric version of the description of the 
{\sl restricted algebra} in [22].
\hfill\qed
\smallskip

\noindent Furthermore, we show that the inequality (9) of 0.2 has a compact counterpart, i.e. when the couple 
$(\cX, \cX_0)$ is replaced by $(X, S)$, where we denote $S\subset X$ a non-singular hypersurface of the projective manifold $X$. The bundle $L\to X$ is assumed to have the properties $\bf(a)-\bf(c)$ above; {\sl in addition,} we assume that we have
$$\alpha\geq \gamma\Theta_h\big(\cal O(S)\big),\leqno(\dagger)$$
where $\gamma$ is a positive real, and $h$ is a non-singular 
metric on the bundle 
$\cal O(S)$ associated to $S$.
 
The hypothesis concerning $\displaystyle \{K_X+S+ {1\over p}L\}$, its corresponding 
minimal metric $\varphi_{\rm min}$ and the metric $\varphi_0$
on $K_X+S+ {1\over p}L_{|S}$ encoded in relations (3)-(8)
are assumed to hold transposed in the actual setting.
In this case, the perfect analogue of (9) is true, as follows: we have
$$\varphi_{\rm min|S}\geq C+ \varphi_0\leqno(11)$$
as it is shown by theorem $\rm B.9$. If we assume that $\nu^j< 1$, then the analogue of the inequality (11) in the present context can be rephrased as follows: {\sl a metric $\varphi_0$ on $\displaystyle K_{X}+ S+ {1\over p}L_{|S}$
is more singular than the restriction of $\varphi_{\rm min}$ to $S$ if and only if it satisfies the relation $\rm (7)$}.
\smallskip
\noindent In section C we prove an extension statement which was used in our previous work [4]: it is a quick consequence of the qualitative version of inequality (12) above (cf. the comments of B.10).

As far as the organization of the present text is concerned, we mention here that the main 
sections $A$ and $B$ can be followed independently (even if they 
share many similar techniques). 
\medskip

\noindent In conclusion, we believe that the {\sl metric point of view} as initiated by J.-P.Demailly in [11] will be extremely useful for further research around the
topics presented in this article.\hfill\qed

\bigskip

\section{ \S A. Proof of theorem 0.1}

\vskip 10pt
\noindent At the beginning of this paragraph we will prove the qualitative part of the theorem 0.1; 
the method we will use it is still the "standard one" borrowed from the articles
in the field quoted above. Nevertheless, there  are quite a few things to be changed 
and therefore we will provide a  complete treatment.

The main technical tool which will be needed is the following effective extension 
theorem. Results of this kind first appeared in [32]; the version which is best adapted for
what we need is taken from [41]. 

\claim A.0 Theorem {\rm ([41])}|Let $\pi:\cX\to \Delta$ 
be a projective family of smooth manifolds. Let $E\to \cX$ be a line bundle,
endowed with a (possibly singular) metric $h$, with semi-positive curvature 
current. If $\displaystyle u\in H^0(\cX_0, K_{\cX_0}+ E)$ is a section of the
adjoint bundle of $E$ restricted to the central fiber,
such that 
$$\int_{\cX_0}|u|^2e^{-\varphi}< \infty$$
then there exist a section $U\in H^0(\cX, K_{\cX}+ E)$, such that 
$\displaystyle U_{\vert \cX_0}= u\wedge d\pi$ and moreover

$$\int_{\cX}|U|^2e^{-\varphi}\leq C_0\int_{\cX_0}|u|^2e^{-\varphi}.$$

\endclaim
\medskip

We recall that the constant $C_0$  above is absolutely universal--in fact, this is the {\sl real strength} of the preceding result.

Anyway, the above result shows that in order to extend some section $u$ of the bundle 
$pK_{{\cal X}_0}+ qL$, it would be enough to
get a metric on  the bundle $(p-1)K_{\cX}+ qL$ such that $u$ is square
integrable with respect to the restriction of this metric to the central fiber. If $p=1$, then $q$ is either equal to 0 or 1 and the metric in question is easy to obtain, since we just take the one we have on the bundle $L$ by hypothesis. In general, the construction of the metric is done by induction, and will be performed in the next two subsections.

\subsection{\S A.1 Choice of the ample line bundle}

\noindent  Let 
$\displaystyle u \in H^0({\cal X}_0, pK_{{\cal X}_0}+ qL))$ be the section we want to extend.
In order to start the inductive procedure which will construct our metric on the bundle
$(p-1)K_{\cX}+ qL$ we first choose an ample line bundle on $\cX$ with the following properties:
\medskip
{\itemindent 8mm
\item {$(A_0)$} For each $\alpha= 0,...,p-q$, the bundle $\alpha K_{\cX}+ A$ is generated by global sections
$(\tau^{(\alpha)}_i)$, where $i= 1,..., Q_\alpha$;
\smallskip
\item {$(A_1)$} For each $\beta= 1,...,q$, the bundle 
$\beta A$ is globally generated by $(s^{(\beta)}_j)$, where $j= 1,..., N_\beta$;
\smallskip
\item{$(A_2)$}  Every section of $\displaystyle {pK_{\cX}+ qL+ (1+q)A}_{\vert {\cX}_0}$
extends to $\cX$;
\smallskip
\item{$(A_3)$} The sheaf $\displaystyle 
{\cal O}(K_X+ L+ A)\otimes {\cal I}(h_{\vert {{\cal X}_0}})$
is generated by its global sections $(s^{(\gamma)})$, for $\gamma= 1,..., M$.

}
\medskip

\noindent Concerning the existence of such a line bundle, see e.g. [40]. 
Remark that $A$ will depend on $(p, q)$
because we impose the extension property $(A_2)$.

\subsection{\S A.2 Inductive procedure}

\vskip 7pt
\noindent We consider a triple of positive integers $(k, \beta, \alpha)$ such that $1\leq \beta\leq q$
and such that $0\leq \alpha\leq p-q$.
In what will follow, we denote by $J$ a {\sl collection} of elements of the set $\{1,...,M\}$
(i.e. we allow repetitions among the elements of $J$) and 
we denote by $s^{(J)}:= \prod _{\rho\in J}s^{(\rho)}$ (we use the notations in the previous paragraph, $A_0-A_3$). The number of the elements of a collection $J$
(including repetitions) will be denoted by $|J|$.

\noindent In order to set-up the inductive procedure, we introduce the next notations.

\noindent 
$\bullet$ If $\beta:= |J|\leq q-1$, then let $u^{(k, J)}_{(j, i)}\in 
H^0\bigl(\cX_0, k(pK_{\cX}+ qL)+ \beta(K_{\cX}+ L)+ (1+q)A_{|\cX_0}\bigr)$
be the section defined by 
$$u^{(k, J)}_{(j, i)}:= u^k\otimes s^{(J)}\otimes s^{(q-\beta)}_j\otimes \tau^{(0)}_i,\leqno (12)$$
where $i= 1,..., Q_0$ and $j= 1,...N_\beta$.

\noindent 
$\bullet$ If $|J|= q$, then let $u^{(k, J, \alpha)}_{(i)}\in 
H^0\bigl(\cX_0, k(pK_{\cX}+ qL)+ q(K_{\cX}+ L)+ \alpha K_\cX+ (1+q)A_{|\cX_0}\bigr)$
be the section defined by 
$$u^{(k, J, \alpha)}_{(i)}:= u^k\otimes s^{(J)}\otimes \tau^{(\alpha)}_i,\leqno(13)$$
where $i= 1,..., Q_\alpha$.

\smallskip

\noindent We formulate the next proposition.
\smallskip
\noindent ${\cP}(k, \beta, \alpha)$:  {\sl Given a triple $(k, \beta, \alpha)$ as above, for any collection $J\subset \{1,...,q\}$ with 
$\vert J\vert = \beta$ we have:

{\itemindent 4mm
\item {$\bullet$} If $1\leq |J|\leq q-1$, then for any $i, j$ as above 
there exists a section
$$U^{(k, J)}_{(j, i)}\in 
H^0\bigl(\cX, k(pK_{\cX}+ qL)+ \beta(K_{\cX}+ L)+ (1+q)A\bigr)$$
whose restriction to the central fiber is equal to
$u^{(k, J)}_{(j, i)}$.

\item {$\bullet$} If $|J|= q$, then for any $i$ there exists a section
$$U^{(k, J, \alpha)}_{(i)}\in 
H^0\bigl(\cX, k(pK_{\cX}+ qL)+ q(K_{\cX}+ L)+ \alpha K_\cX+ (1+q)A\bigr)$$
whose restriction to the central fiber is equal to $u^{(k, J, \alpha)}_{(i)}$.
\hfill\qed

}}

\bigskip
\noindent The core of the proof lies in the next statement (see [35]).
\claim Lemma A.2.1|The proposition ${\cP}(k, \beta, \alpha)$ is true for any $k\in \bZ_+$, any $1\leq \beta\leq q$ and any $0\leq \alpha\leq p-q$.
Moreover, the effective version of ${\cP}(k, \beta, \alpha)$ holds. There exists a constant $C> 0$  independent of $k$ such that if we denote by $h^{(k, \beta)}$ 
--respectively $h^{(k, q, \alpha)}$-- the algebraic metric on the bundle
$k(pK_{\cX}+ qL)+ \beta(K_{\cX}+ L)+ (1+q)A$ --respectively on 
$k(pK_{\cX}+ qL)+ q(K_{\cX}+ L)+ \alpha K_{\cX}+(1+q)A$--
induced by the family of sections  
$\big(U^{(k, J)}_{(j, i)}\bigr)_{ i, j, |J|=\beta\leq q-1}$ --resp. by $\big(U^{(k, J,\alpha)}_{(i)}\big)_{i, |J|= q}$-- then we have
$$\max\Big\{\int_\cX e^{\varphi^{(k, \beta)}- \varphi^{(k, \beta-1)}-\varphi_L},
\int_\cX e^{\varphi^{(k, 1)}- \varphi^{(k-1, q, p-q)}} \Big\}\leq C\leqno(14)$$
and 
$$\max\Big\{\int_\cX e^{\varphi^{(k, q, \alpha)}- \varphi^{(k, q, \alpha-1)}}, 
\int_\cX e^{\varphi^{(k, q, 0)}- \varphi^{(k-1, q-1)}-\varphi_L}
 \Big\}\leq C.\leqno(15)$$
\endclaim

\proof.  Even if the formulation of the above lemma is somehow complicated, the inductive procedure is quite natural and it will be performed as follows: 
we first show that ${\cP}(1, 0, 0)$ is valid, and then we prove that the next sequence of implications holds true
$${\cP}(1, 0, 0)\to {\cP}(1, 1, 0)\to {\cP}(1, 2, 0)\to...\to {\cP}(1, q, 0)\to 
$$
$$\to{\cP}(1, q, 1)\to{\cP}(1, q, 2)\to ...\to {\cP}(1, q, p-q)\to{\cP}(2, 1, 0)\to {\cP}(2, 2, 0)\to...$$
Notice that we allow $\beta$ to be equal to $0$ only for the first term in the previous sequence of implications: the reason is that $\varphi^{(k-1, q, p-q)}= 
\varphi^{(k, 0)}$.
We remark that even if the proposition ${\cP}(k, \beta, \alpha)$ is purely qualitative, the procedure we describe next will
produce the uniform constant $"C``$ in the statement above as well.
\medskip

To check the first proposition ${\cP}(1,0,0)$ is fairly easy: it is just the fact that $A$ is positive enough
to satisfy the property $(A_2)$; this allows the extension of the sections $u\otimes s^{q}_j\otimes \tau^{(0)}_i$
for each $j= 1,..., N_q$ and $i=1,..., Q_0$.

Assume now that for some indexes $(k, \beta, \alpha)$ the property ${\cP}(k, \beta, \alpha)$ has been established. Then we have to distinguish 
between several cases.
\smallskip
$\bullet$  We first consider the case $\alpha= 0$ and $1\leq \beta\leq q-1$. By property ${\cP}(k, \beta, 0)$ we deduce that for each 
indexes 
$(i, j, J)$ such that $1\leq j\leq N_{q-\beta}, 1\leq i\leq Q_0, |J|= \beta$, the section 
$$u^{(k, J)}_{j, i}\in 
H^0\bigl(\cX_0, k(pK_{\cX_0}+ qL)+ \beta(K_{\cX_0}+ L)+ (1+q)A\bigr)$$ 
admits an extension
$$U^{(k, J)}_{(j, i)}\in 
H^0\bigl(\cX, k(pK_{\cX}+ qL)+ \beta(K_{\cX}+ L)+ (1+q)A\bigr).$$

\noindent Next we use the family of sections $\displaystyle \big(U^{(k, J)}_{j,i})_{1\leq j\leq N_{q-\beta}, 1\leq i\leq Q_0, |J|= \beta}$ to construct a metric 
$h^{(k,\beta)}$ on the bundle 
$$k(pK_{\cX}+ qL)+ \beta(K_{\cX}+ L)+ (1+q)A;$$ 
it will be singular in general, but its singularities over the central fiber are perfectly understood.  

\noindent For each collection of integers $K$ such that $|K|= \beta+1$ and for each integers 
$i,j$ let us consider the section 
$$u^{(k, K)}_{j, i}\in 
H^0\bigl(\cX_0, k(pK_{\cX_0}+ qL)+ (1+ \beta)(K_{\cX_0}+ L)+ (1+q)A\bigr);$$
we intend to extend it in an effective manner by Ohsawa-Takegoshi theorem. To this end, we decompose the bundle above as follows
$$k(pK_{\cX}+ qL)+ (1+ \beta)(K_{\cX}+ L)+ (1+q)A= K_{\cX}+ k(pK_{\cX_0}+ qL)+ \beta(K_{\cX_0}+ L)+ (1+q)A+ L$$
and remark that in this way it becomes the adjoint bundle of 
$$E:= k(pK_{\cX}+ qL)+ \beta(K_{\cX}+ L)+ (1+q)A+ L.$$
Now the bundle $E$ can be endowed with the metric $\displaystyle h^{(k,\beta)}\otimes h_L$; it is 
semi-positively curved, and we check now the integrability of the 
section we want to extend with respect to it. We have the next relations.
$$\eqalign{
&I: = \int_{\cX_0}
|u^{(k, K)}_{j, i}|^2e^{-\varphi^{(k,\beta)}-\varphi_L}=
\int_{\cX_0}{{|u^{(k, K)}_{j, i}|^2}\over {
 \sum_{l, m} ^{|J|= \beta} {|U^{(k, J)}_{l, m}|^2}}}e^{-\varphi_L}= \cr
&= \int_{\cX_0}{{|u^{(k, K)}_{j, i}|^2}\over {
 \sum_{l, m} ^{|J|= \beta} {|u^{(k, J)}_{l, m}|^2}}}e^{-\varphi_L}  \leq C\int_{\cX_0}{{\big(\sum_{\gamma}|s^{(\gamma)}|^2\big)^{\beta+1}}\over
{\big(\sum_{\gamma}|s^{(\gamma)}|^2\big)^{\beta}}}e^{-\varphi_L-\varphi_A}\leq \cr
&\leq C.\cr 
}$$

The second equality holds because of the definition 
of the metric $h^{(k, \beta)}$; the third one
is given by the extension property ${\cP}(k, \beta, 0)$. The fourth inequality is obtained 
by simplification of the 
the common factor $u^k$, and the fact that all indexes $J$ such that $|J|= \beta$ appears in the expression of the denominator. We also use the fact that the sections $(\tau^{(0)}_i)$ do not have common zeroes.
Finally the last inequality follows from the fact that the sections $s^{(\gamma)}$ belong to the multiplier ideal of the restriction of the 
metric $h_L$ to the central fiber.

 The constant $"C``$ in the last line only depends on the auxiliary sections $(s^{(\gamma)}, s^{(m)}_j)$ and thus they are
uniform with respect to $k$; also $\varphi_A$ is just any smooth metric on $A$.
 
\noindent  Thus the requirements of the extension theorem 2.1 are satisfied, 
and therefore for each indexes $(K, i, j)$ we obtain 
$$U^{(k, K)}_{(j,i)}\in H^0\bigl(\cX, k(pK_{\cX}+ qL)+ (1+ \beta)(K_{\cX}+ L)+ (1+q)A$$ 
such that:

\item {(i)} ${U^{(k, K)}_ {(j,i)}}_{|\cX_0}= u^{(k, K)}_ {(j,i)}$;
\smallskip
\item{(ii)} We have
$$\int_\cX e^{\varphi^{(k, \beta)}- \varphi^{(k, \beta-1)}-\varphi_L}\leq C$$
for some constant $C$ which is a fixed multiple of one obtained a few lines above. Indeed, all we have to do is to add up
the several estimates obtained above, and remark that the number of the terms is 
bounded uniformly with respect to $k$.

\noindent Therefore, the first case is completely settled.\hfill \qed
\vskip 5pt
\smallskip
\noindent $\bullet$ We analyze here the second case, namely $\alpha= 0$ and 
$\beta= q$; the arguments are quite similar to the previous case. Since we admit the validity of $\cP(k, q, 0)$, we have the
family of sections 

$$U^{(k, J)}_{(i)}\in H^0\bigl(\cX, k(pK_{\cX}+ qL)+ q(K_{\cX}+ L)+ (1+q)A\bigr)$$ 
such that
$${U^{(k, J)}_{(i)|\cX_0}}{}= u^{(k, J)}_{(i)};$$
as before, we can use them to define a metric $h^{(k, q)}$ on the bundle 
$$k(pK_{\cX}+ qL)+ q(K_{\cX}+ L)+ (1+q)A.$$ 
We have to extend each member of the family of sections 
$$u^{(k, K, 1)}_{(i)}\in H^0\bigl(\cX_0, 
k(pK_{\cX}+ qL)+ q(K_{\cX}+ L)+ K_{\cX}+ (1+q)A_{|\cX_0}\bigr),$$
where $|K|= q$. To this end
we will use again the Ohsawa-Takegoshi theorem; we can write
$$k(pK_{\cX}+ qL)+ q(K_{\cX}+ L)+ K_{\cX}+ (1+q)A= K_{\cX}+ k(pK_{\cX}+ qL)+ q(K_{\cX}+ L)+ (1+q)A$$
and remark that in this way it become the adjoint bundle of 
$$E:= k(pK_{\cX}+ qL)+ q(K_{\cX}+ L)+ (1+q)A.$$
The bundle $E$ can be endowed with the metric $\displaystyle h^{(k, q)}$; it is 
semi-positively curved, and we check now the integrability of the 
section above.

$$\eqalign{
&I: = \int_{\cX_0}
|u^{(k, K, 1)}_{i}|^2e^{-\varphi^{(k, q)}}=
\int_{\cX_0}{{|u^{(k, K, 1)}_{i}|^2}\over {
 \sum_{l, m} ^{|J|= q} {|U^{(k, J)}_{l, m}|^2}}}= \cr
&= \int_{\cX_0}{{|u^{(k, K, 1)}_{i}|^2}\over {
 \sum_{l, m} ^{|J|= q} {|u^{(k, J)}_{l, m}|^2}}} \leq C\int_{\cX_0}{{\big(\sum_{\gamma}|s^{(\gamma)}|^2\big)^{q}}\over
{\big(\sum_{\gamma}|s^{(\gamma)}|^2\big)^{q}}}dV\leq \cr
&\leq C.\cr 
}$$
Thus, the second case is completely solved.\hfill \qed

\bigskip
\noindent $\bullet$ The remaining cases we have to consider are $(1\leq \alpha\leq p-q-1,  \beta= q)$ and respectively $(\alpha= p-q, \beta= q)$. We only give the arguments for the latter (and we leave the former to the interested reader).  
The implication we have to prove is 
$${\cP}(k, q, p-q)\to{\cP}(k+1, 1, 0)$$
Since the proposition ${\cP}(k, q, p-q)$ is valid, we have the family of sections 
$$U^{(k, J, p-q)}_{(i)}\in 
H^0\bigl(\cX, (k+ 1)(pK_{\cX}+ qL)+ (1+q)A\bigr)$$ 
such that 
$$U^{(k, J, p-q)}_{{(i)}|\cX_0}= 
u^{(k, J, p-q)}_{(i)}$$
Let $h^{(k, q, p-q)}$ be the algebraic metric given by the sections 
$\displaystyle (U^{(k, J, p-q)}_{(i)})$ above, where $|J|= q$ and 
$1\leq i\leq Q_{p-q}$. Consider the section
$$u^{(k+1, K)}_{(j, i)}\in 
H^0\bigl(\cX_0, (k+1)(pK_{\cX_0}+ qL)+ K_{\cX_0}+ L+ (1+q)A\bigr)$$ 
(where $|K|= 1$). We check now its integrability with respect to the metric $h^{(k, q, p-q)}$
twisted with the metric of $L$; in the forthcoming computations we skip some trivial steps which are direct consequences of the definition of the corresponding objects.

$$\eqalign{
I: & = \int_{\cX_0}|u^{(k+1, K)}_{(j, i)}|^2e^{-\varphi^{(k,q)}-\varphi_L)}\leq  C\int_{\cX_0}
{{|u\otimes s^{(K)}\otimes s^{(q-1)}_j\otimes\tau^{(0)}_i|^2}\over {
 \sum_{l, |J|= q} |s^{(J)}\otimes\tau^{(p-q)}_l|^2}}\exp(-\varphi_L)
\leq \cr
& \leq C.\cr 
}$$

Remark that the last integral converge precisely because of the hypothesis
$$u\in \displaystyle H^0\Bigl(\cX_0, pK_{\cX_0}+ qL)\otimes 
\ol{\cI}
(h_{L\vert \cX_0})^{(q)}\Bigr)$$ 
 and this ends the proof of the lemma.\hfill\qed
\vskip 10pt

\noindent The estimates (14) and (15) of lemma A.2.1 show that we can consider the limit metric
$$h^{(\infty)}:= \lim_k {h^{(k, 1)}}^{1\over k}$$
of the bundle $pK_{\cX}+ qL)$, which in addition has the next properties:

\item{(a)} The curvature current of $h^{(\infty)}$ is positive;
\smallskip
\item {(b)} The restriction of the metric $h^{(\infty)}$ to the central fiber
is well-defined, and we have 
$\displaystyle \sup _{\cX_0}|u|_{h^{(\infty)}}< \infty$.

\noindent For the existence of the limit and the verification of the above relations 
we refer e.g. to [15].
\medskip
Now a last application of Ohsawa-Takegoshi extension result will show that the section 
$u $ extend over the whole family. Indeed, we have
$$pK_{\cX}+ qL= K_{\cX}+ {{p- 1}\over {p}}(pK_{\cX}+q L)+ {q\over p}L$$
and we endow the bundle $\displaystyle {{p- 1}\over {p}}(pK_{\cX}+q L)+ {q\over p}L$ with the metric 
$\bigl(h^{(\infty)}\bigr)^{1-1/p}\otimes h_L^{q\over p}$. We have

$$\int_{\cX_0}
|u|^2e^{-{p-1\over p}\varphi^{{(\infty)}} -{q\over p}\varphi_L}
\leq C_0\int_{\cX_0}
|u|^{2/p} e^{-{q\over p}\varphi_L}<\infty$$
where the first inequality is given by the property $(b)$ above, and for the last one we use 
the fact that the coefficients of the section $u$ belong to the ideal $\ol{\cal I}
(h_{L\vert \cX_0})^q$, together with the H\"older inequality. Again, we see that the condition
$p\geq q$ is crucial.

This finishes the proof of the theorem 0.1, modulo the integrability of the appropriate 
root of the extension. To clear this last point, we first show that 
$$\int_{\cX}e^{{\varphi^{(\infty)}- q\varphi_L\over p}}< \infty.\leqno (16)$$
The relation above is obtained as follows: we multiply the inequalities (14) and (15)
for successive parameters, and we use 
H\"older inequality. We infer the existence of a positive constant $C$ such that 
 $$\int_{\cX}e^{{{\varphi^{(k+1, 1)}- \varphi^{(1, 1)}}
 \over {kp}}- {q\over p}\varphi_L}d\lambda\leq C \leqno (17)$$
 for any $k\geq 1$. 
 
 Elementary properties of
 plurisubharmonic functions show that the sequence 
 $${1\over k}\varphi^{(k, 1)}$$ 
 converges 
 a.e. and in $L^1$ to the metric $\varphi^{(\infty)}$ (up to the choice of a sub-sequence). 
 By Lebesgue's dominated convergence we can take the limit as $k\to\infty$ in (16) and obtain (17).
 
 Now remember that the extension $U$ of our section $u$ satisfies the following 
 $L^2$ estimate
 $$\int_{\cX}|U|^2e^{-{p-1\over p}\varphi^{{(\infty)}} -{q\over p}\varphi_L}
< \infty.$$ 
Then we have
$$\leqno(18)
\eqalign{
\int_{\cX} & |U|^{2\over p}e^{-{q\over p}\varphi_L}=  
\int_{\cX} |U|^{2\over p}e^{-{{p-1}\over {p^2}}\varphi^{{(\infty)}}-{{q}\over {p^2}}\varphi_L}e^{{{p-1}\over {p^2}}\varphi^{{(\infty)}}+
{{(1-p)q}\over {p^2}}\varphi_L}\leq \cr
 \Bigl(\int_{\cX} & |U|^2e^{-{p-1\over p}\varphi^{{(\infty)}} -{q\over p}\varphi_L}
\Bigr)^{{{1}\over {p}}}\Bigl(\int_{\cX}
e^{{\varphi^{(\infty)}- q\varphi_L\over p}}\Bigr)^{1-{{1}\over {p}}}< \infty
\cr
}$$

\noindent The first part of the proof of 0.1 is now complete.
\hfill \qed 

\medskip
\noindent In the last part of this section, we establish the quantitative part of the theorem 0.5, namely the existence of 
a section 
$$U\in \displaystyle  
H^0\bigl(\cX, pK_{\cX}+ qL\bigr)$$ 
such that: 

{\itemindent 5mm
\item {(i)} Over the central fiber we have $U_{|\cX_0}= u\wedge d\pi^{\otimes p}$;
\smallskip
\item {(ii)} The next
$L^{1/p}$ integrability condition holds
$$\int_{\cX}|U|^{2\over p}e^{-{q\over p}\varphi_L}\leq C_0 \int_{\cX_0}[|u|^{2\over p}e^{-{q\over p}\varphi_L}.
$$}


\proof \hskip 1mm {\sl of (ii)}.  We will use basically the same arguments as in the proof of the 
$\displaystyle L^{2\over m}$ extension theorem in the paper [4].

In the first place we observe that the space of all the possible extensions of $u$ with integrable $L^{2/p}$ semi-norm 
is non-empty, thanks to (18)--this is the crucial point!
Next we define $U$ to be an extension of $u$ which minimize the previous semi-norm;
with this choice we show now that the estimate required in the theorem above is satisfied. 

Indeed, let us consider  the bundle 
$$pK_{\cX}+ qL= K_{\cX}+ {{p-1}\over {p}}(pK_{\cX}+ qL)+ {q\over p}L;$$
it is the adjoint bundle of $\displaystyle {{p-1}\over {p}}(pK_{\cX}+ qL)+ {q\over p}L$, and we can endow the latter with
the metric induced by the section $U$ raised to the power $1-1/p$, twisted with the metric of $q/pL$.
This metric has semi-positive curvature and can be restricted to the central fiber, 
as it is the case for the metric of $L$, and the section $U$ is not identically zero on $\cX_0$. 

The section $u\in \displaystyle H^0\bigl(\cX_0, pK_{\cX_0}+ qL\bigr)$ is square integrable with respect to
the previous metric, because the integrability condition reads as
$$\int_{\cX_0} {|u|^2\over |u|^{2{p-1\over p}}}e^{-{q\over p}\varphi_L}= \int_{\cX_0}
|u|^{2\over p}e^{-{q\over p}\varphi_L}< \infty . $$
We use here the fact that $U$ is an extension of $u$, as well as the hypothesis that $u$ belongs to the appropriate 
power of the multiplier ideal sheaf, which implies in particular that the last integral above is finite. 

Thus the Ohsawa-Takegoshi theorem shows the existence of some extension 
$$U_1\in \displaystyle  H^0\bigl(\cX, pK_{\cX}+ qL\bigr)$$ 
of our section $u$ such that 
$$\int_{\cX} {|U_1|^2\over |U|^{2{p-1\over p}}}e^{-{q\over p}\varphi_L}\leq C_0
\int_{\cX_0}|u|^{2\over p}e^{-{q\over p}\varphi_L}.$$
But then we are done, since we necessarily have  
$$\int_{\cX}|U|^{2\over p}e^{-{q\over p}\varphi_L}\leq 
\int_{\cX} {|U_1|^2\over |U|^{2{p-1\over p}}}e^{-{q\over p}\varphi_L}$$
because if not the minimality property of the section $U$ will be violated: the argument is as follows. We assume that the inequality above does not hold; then we have

$$\eqalign{
& \int_{\cX}|U_1|^{2\over p}e^{-{q\over p}\varphi_L}= 
\int_X {{|U_1|^{2\over p} e^{-{{q\varphi_L }\over {p^2}}}\over {|U|^{2{p-1\over p^2}} }}
|U|^{2{p-1\over p^2}}} 
e^{-{q\over p^2}(p-1)\varphi_L} \leq \cr
 & \Big(\int_{X} {|U_1|^2\over |U|^{2{p-1\over p}}}e^{-{q\over p}\varphi_L}\Big)^{1\over p}
\Big(\int_{X} |U|^{2\over p}e^{-{q\over p}\varphi_L}\Big)^{p-1\over p} \cr
 & < \int_{X} |U|^{2\over p}e^{-{q\over p}\varphi_L}.\cr
}$$
The contradiction we have just obtained shows our result 0.1 is completely proved.
\hfill \qed


\section {\S B. Canonical metrics and their restriction properties}

\smallskip
\noindent In this section we prove theorem 0.2 and we derive some of its consequences.

Let $\pi:\cX\to \bD$ be a proper, surjective map, where $\bD$ is the unit disk.
We assume that the central fiber $\cX_0= \pi^{-1}(0)$ is non-singular, and let $L\to \cX$
be a hermitian line bundle such that we have
$$p([\Delta]+ \alpha)\in c_1(L);\leqno(19)$$
we assume that $\Delta= \sum_{j\in J}\nu^jY_j$ and the metric 
$e^{-\varphi_0}$ on $K_{\cX}+ {1\over p}L_{|\cX_0}$ are 
satisfying the properties $\bf(a)-(c)$, respectively (7), (8) in the introduction. 
The conventions in the introduction are in force during all of the present paragraph.

\smallskip

\noindent Following Ein-Popa's elegant approach in [16], for each $s=1,..., p$ we define the set
$$J_s:= \{j\in J^\prime : p\nu^j\geq s\}.$$
Then we
can write
$$p\sum_{j\in J^\prime}\nu^jY_j= \sum_{s=1}^{p}\sum_{j\in J_s}Y_j,\leqno(20)$$
and the relation (20) induces a decomposition
$$L= L_1+...+ L_{p-1}+ L_p$$
where $L_s\equiv \sum_{j\in J_s}Y_j$ for each $s= 1,..., p-1$, 
and such that $L_p$ admits a metric whose curvature form equals
$p(\alpha+ \sum _{j\in J\setminus J^\prime}\nu^j[Y_j])+ \sum_{j\in J_p}
[Y_j]$.

Let $k\in \bZ_+$ and $r\in \{0, 1, ..., p-1\}$; we introduce the notation 
$$L^{(r)}:= rK_{\cX}+ L_1+...+ L_r\leqno(21)$$
together with the convention that $L^{(0)}$ is the trivial bundle. 
By [40], there exists an ample line bundle $A$ on $\cX$  having the following uniform global generation property: {\sl for any positively curved hermitian bundle $(F, h_F)$ on the central fiber $\cX_0$, the sheaf 
$$\cO\big((K_{\cX_0}+ F+ L^{(r)}+ A-\sum_{j\in J}Y_j)\otimes \cI(h_F)\big)\leqno(22)$$
is generated by its global sections, for any $r= 0,..., p-1$.} 
We also assume that $A$ is ample enough, so that the bundles
$L^{(r)}+ A$, their adjoints $K_{\cX}+ L^{(r)}+ A$ as well as $L^{(r)}+ A-\sum_{j\in J}Y_j$ are very ample, for $r= 0,..., p$.
\hfill\qed
\smallskip  
\noindent We introduce next the main technical tool which will lead us to 0.2. 

For each $(k, r)$ within the range prescribed above we will briefly recall the construction of the $kp$-Bergman metric on the bundle 
$$k(pK_{\cX/\bD}+ L)+ L^{(r)}+ A$$
where we denote by $K_{\cX/\bD}:= K_{\cX} - p^*K_{\bD}$ the relative canonical bundle of the map $\pi: \cX\to \bD$. The existence of this metric, together with its main features which are included in the next statement 
is crucial for the proof of 0.2.
We refer to the articles [3], [4] for details and proofs (see also [47], [48] and the references therein for related results).

Let $\bD^\prime\subset \bD$ be a Zariski open set, such that for each 
$t\in \bD^\prime$, each section of the bundle
$$k(pK_{\cX}+ L) + L^{(r)}+ A_{|\cX_t}$$
extends locally near $t$, for all $r= 0,... p-1$. 
The result proved in [4] states as follows.

\claim B.1 Theorem([4])|There exists a positively curved metric 
$h^{(kp+ r)}_{\cX/\bD}$ on the bundle $$k(pK_{\cX/\bD}+ L)+ L^{(r)}+ A\leqno (23)$$ such that 

\item {\rm (a)} For any $t\in \bD^\prime$, the restriction of the dual metric $h^{(kp+r)\star}_{\cX/\bD}$ to $\cX_t$ is defined by 
$$|\xi|:= \sup_{u\in B_t^{kp}(1)}|\xi (\wt u_x)|$$
where $\xi$ is a vector in the dual bundle fiber $-k(pK_{\cX/\bD}+ L)- L^{(r)}- A_{|\cX_t, x}$. We denote by $B_t^{kp}(1)$ the set of all holomorphic sections $u$ of the bundle 
$\rm (23)$
restricted to ${\cX_t}$ satisfying 
$$\int_{\cX_t}|u|^{2/kp}\exp \Big(-{{\varphi_{r, A}}\over {kp}}-{{\varphi_{L}}\over {p}}\Big)d\lambda\leq 1,$$
and we denote by $\wt u:= u\wedge d\pi^{\otimes kp}$. The metric $\varphi_{r, A}$ is non-singular, positively curved on $L^{(r)}+ A$, and $\varphi_L$ is induced by {\rm (19)}.
\smallskip
\item {\rm (b)} For each compact set $K\subset \Delta$ there exists 
a constant $C_K> 0$ uniform with respect to $k$, such that the local weights $\varphi^{(kp+r)}_{\cX/\bD}$ of the metric 
$h^{(kp+ r)}_{\cX/\Delta}$ are bounded from above by $kC_K$ on every co-ordinate set contained in $\pi^{-1}(K)$. \hfill\qed 

\endclaim

\noindent An important observation is that the metric constructed above {\sl is not} explicitly described
on the set $\bD\setminus \bD^\prime$, so a priori we don't know the size of its singularities over that set. However, as we have remarked in our previous article [4] 
the ``extendable sections" of the restriction
$$k(pK_{\cX/\bD}+ L)+ L^{(r)}+ A_{|\cX_t}$$ 
provides us with a {\sl lower bound} for the weights of $h^{(kp+ r)}_{\cX/\Delta}$, even if
$t\in \bD\setminus \bD^\prime$. Let us explain this next.

The main claim is the following. Let $\mu> 0$ be a real number such that the disk centered at zero with radius $\mu$ does not contain any critical value of $\pi$, and let 
$\tau\in \bD$ such that $|\tau|< \mu$.
We consider 
a holomorphic section $U$ of the bundle 
$k(pK_{\cX/\bD}+ L)+ L^{(r)}+ A$ over the whole family $\cX$, whose global $L^{2/kp}$ norm is finite;
then (modulo an abuse of notation) we have 
$${{|U(x)|^2e^{-\varphi^{(kp+r)}_{\cX/\bD}(x)}}\over {\Big(\int_{\cX_\tau}|U|^{2/kp}\exp \big(-{{\varphi_{r, A}}\over {kp}}-{{\varphi_{L}}\over {p}}\big)d\lambda\Big)^{kp}}}\leq 1
\leqno (24)$$
where $x\in \cX_\tau$ is an arbitrary point.

Indeed, if $\tau\in \bD^\prime$, then the above claim is a consequence of the definition. If not, then use a limit argument--since the weights 
$\varphi^{(kp+ r)}_{\cX/\Delta}$ are {\sl upper semi-continuous}, and since the singularities of
$1/p\varphi_L$ are mild enough (see [4]). \hfill\qed

\medskip 
\noindent We come back now to the metric $\varphi_0$ given by hypothesis and we use it to define the space  
$$V_{k, r}:= {\rm H}^0\big(\cX_0,  \big(kpK_{\cX_0}+ kL+ L^{(r)}+ A_{|\cX_0}\big)\otimes \cI(\psi_{k, r})\big)\leqno(25)$$
where 
$$\psi_{k, r}:= (kp-1)\varphi_0+ \sum_{j\in J}(1+ \nu^j)\log |f_{Y_j}|^2+ \wt \varphi_{r, A}.\leqno(26)$$ 
In the expression of the metric above, we denote by 
$\wt \varphi_{r, A}$ a non-singular, positively curved metric on the bundle
$L^{(r)}+ A-\sum_{j\in J}Y_j$. We remark that we introduce an additional singularity 
$(1+ \nu^j)$ instead of $\nu^j$
in the expression of the metric $\psi_{k,r}$; it will be useful during the proof of lemma B.3.

Anyway, the 
set $V_{k, r}$ is in fact a Hilbert space, whose inner product is given by the formula
$$ \langle\!\langle u, v\rangle\!\rangle:= \int _{\cX _0} \langle u, v\rangle 
e^{-\psi_{k, r}}.\leqno(27)$$
\medskip

\noindent We consider an orthonormal basis $(u^{(kp+r)}_j)$ of $V_{k, r}$ and we prove next the assertions B.2-B.5, that together will prove theorem 0.2. The approach presented here has many similarities with and generalizes the one in [15], [47], [48].
\claim B.2 Lemma|There exists a constant $C$ independent of $k, j$, such that
$$\int_{\cX_0}
|u^{(kp+r)}_j|^{2/kp}\exp \big(-{{1}\over {p}}\varphi_{L}-{{1}\over {kp}}\varphi_{r, A}\big)\leq C
$$
for all $k\gg 0$.
\endclaim

\proof. This is a consequence of H\"older inequality, as follows.

$$\eqalign{
 \int_{\cX_0}& 
|u^{(kp+r)}_j|^{2/kp} e^{-{{1}\over {p}}\varphi_{L}-{{1}\over {kp}}\varphi_{r, A}} =  \cr
= \int_{\cX_0} &
|u^{(kp+r)}_j|^{2\over kp} e^{-{kp-1\over kp}\varphi_0-{{1}\over {kp^2}}\varphi_L- {{1}\over {kp}}\varphi_{r, A}} 
e^{{kp-1\over kp}\varphi_0-{kp-1\over kp^2}\varphi_{L}} \cr
\leq &\Vert u^{(kp+r)}_j\Vert ^{1\over kp} \Big(\int_{\cX_0} e^{\varphi_0-{1\over p}\varphi_{L}}\Big)^{kp-1\over kp}\leq C.\cr
}$$
The last inequality is valid because of the integrability assumption (9) concerning the metric $\varphi_L$. \hfill\qed
\medskip
We introduce the set 
$$J^1:= \{j\in J : \nu^j= 1\}$$
and we show next that each element of $V_{k, r}$ admits an extension to $\cX$
which vanishes along the divisor $\displaystyle \sum_{j\in J^1}Y_j$. This will be crucial for the study of $\varphi^{(kp+r)}_{\cX/\bD}$, given the inequality (24). Most of the ``extension" arguments provided for the following lemma has been invented in [40]; to our knowledge, their relevance in the actual context first appeared in [47]. 
\claim B.3 Lemma|For each $k, r$ and $j$ there exists a section 
$$U^{(kp+r)}_j\in H^0\Big(\cX,  kpK_{\cX/\bD}+ kL+ L^{(r)}+ A\Big)$$
whose restriction to $\cX_0$ is equal to $u^{(kp+r)}_j$, and such that its zero divisor contains
$\displaystyle \sum_{j\in J^1}Y_j$. \hfill\qed
\endclaim

\proof. 
We will use induction on $kp+ r$; if $k=1$ and $r= 0$, then 
the extension of the sections $u^{(p)}_j$ is a consequence of the ampleness of $A$, together with the fact that the $L^2$ condition in the definition of the space $V_{1, 0}$ 
force the vanishing of $u^{(p)}_j$ along $\displaystyle \sum_{j\in J}Y_j$. 

\noindent Therefore, we assume that the extension of the sections 
$u^{(kp+r)}_j$ to $\cX$ with the vanishing properties required by B.3 has been already shown to exist, and let 
$u^{(kp+r+1)}_i$ be an element of the basis of $V_{k, r+1}$.

\noindent $\bullet$ If $r\leq p-2$, then we intend to use the global generation property of the bundle $A$ (see (22)) where the data is 
$$F:= (kp-1)\big(K_{\cX_{|0}}+ {1\over p}L\big)+ {1\over p}L+ \sum_{j\in J}Y_j$$
and $\displaystyle \varphi_F:= (kp-1)\varphi_0+ \sum_{j\in J}(1+\nu^j)\log |f_{Y_j}|^2$.
Since the section $u^{(kp+r+1)}_i$ belongs to the ideal associated to the metric $\varphi_F$, 
we
have the {\sl pointwise} inequality
$$|u^{(kp+r+1)}_i|^2\leq C\sum_j |u^{(kp+r)}_j|^2\leqno(28)$$
by global generation property (22), where the norms are computed with respect to some non-singular metric on the corresponding bundle. 

We write the bundle $kpK_{\cX/\bD}+ kL+ L^{(r+1)}+ A$ in adjoint form, as follows
$$kpK_{\cX/\bD}+ kL+ L^{(r+1)}+ A= K_{\cX/\bD}+ L_{r+1}+ kpK_{\cX/\bD}+ kL+ 
L^{(r)}+ A$$
By Ohsawa-Takegoshi theorem A.0, in order to extend the section 
$u^{(kp+r+1)}_i$, it is enough to endow the bundle
$$E:= L_{r+1}+ kpK_{\cX/\bD}+ kL+ 
L^{(r)}+ A$$
with a semi-positively curved metric, such that the $L^2$ norm of $u^{(kp+r+1)}_i$
with respect to it is finite.
We denote by $\varphi_{L_{r+1}}$ the singular metric on $L_{r+1}$, whose curvature current 
is equal to $\sum_{j\in J_{r+1}}[Y_j]$, and by $\wt\varphi_{L_{r+1}}$ a non-singular metric on this bundle -for which we cannot impose any curvature requirements. We also define the metric $h^{(kp+ r)}$ on $kpK_{\cX/\bD}+ kL+ L^{(r)}+ A$ induced by the family of sections 
$\big(U^{(kp+r)}_j\big)_j$.

For any parameters $\delta, \varepsilon, \tau\in \bR_+$ we define the next 
metric on the bundle $E$
$$\varphi_E:= (1-\delta)\varphi_{L_{r+1}}+\delta \wt\varphi_{L_{r+1}}+ 
(1-\varepsilon)\varphi^{(kp+ r)}+ \varepsilon \big((kp-1)\varphi_{\rm min, \tau}+ 
{1\over p}\varphi_L+ \wh\varphi_{r, A}\big)\leqno(29)$$
where $\wh\varphi_{r, A}$ is a positively curved 
non-singular metric on the bundle $K_{\cX}+ L^{(r)}+ A$ and 
$\varphi_{\rm min, \tau}$ is the metric induced on $K_{\cX}+ 1/pL$ by the metric with minimal singularities on $K_{\cX}+ 1/pL+ \tau A$. We note that its curvature form is greater than $-\tau\omega_A$, and that its restriction to 
$\cX_0$ has the expression in (6).

We remark that $(E, e^{-\varphi_E})$ is positively curved, provided that 
$\varepsilon\gg \delta$, and that $(kp-1)\tau\ll 1$.
We still have to check that the following integral is convergent 
$$\int_{\cX_0}|u^{(kp+r+1)}_i|^2e^{-\varphi_E}< \infty.\leqno(30)$$
From relation (28), we see that the above $L^2$ condition will be satisfied 
if we can show that
$$\int_{\cX_0}|u^{(kp+r+1)}_i|^{2\varepsilon}e^{-(1-\delta)\varphi_{L_{r+1}}-\varepsilon
(kp-1)\varphi_{\rm min, \tau }-{\varepsilon\over p}\varphi_L}< \infty\leqno(31)$$
(we ignore the non-singular weights in the expression of the $\varphi_E$). 

In order to establish the relation (31), we recall that by hypothesis we have
$$\varphi_0\leq \sum_{j\in J^\prime}\rho^j_{\rm min, \infty}\log|f_j|^2+ 
\sum_{j\in J\setminus J^\prime}\nu^j\log|f_j|^2.$$
Since the section
$u^{(kp+r+1)}_i$ belongs to the space $V_{k, r+1}$ we infer that the divisor
$$\sum_{j\in J^\prime}\big([(kp-1)\rho^j_{\rm min, \infty}+ \nu^j]+ 1\big)Y_j+
\sum_{j\in J\setminus J^\prime}(kp\nu^j+ 1)Y_j$$
is {\sl smaller than} its zero divisor. 
As recalled in the introduction, for any $\tau> 0$ we have $\rho^j_{\rm min, \infty}\geq \rho^j_{\rm min, \tau}$, and therefore the integral
(31) is dominated by the quantity
$$\int_{\cX_0}{e^{-(1-\delta)\varphi_{L_{r+1}}-\varepsilon
(kp-1)\varphi_{\Lambda_{\tau 0}}}d\lambda\over 
\prod_{j\in J\setminus J^\prime}
|f_j|^{2\varepsilon\big(kp(\rho^j_{\rm min, \infty}-\nu^j)- \rho^j_{\rm min, \infty}-1\big)}}.$$
Indeed, the lower bound of the vanishing of the section $u^{(kp+r+1)}_i$ as 
explained before is big enough in order to compensate the singularities
$$\varepsilon\sum_{j\in J^\prime}\big((kp- 1)\rho^j_{\rm min, \infty}+ \nu^j\big)\log|f_{Y_j}|^2$$
arising from the restriction of $\varphi_{\rm min, \tau}$ and $\displaystyle{1\over p}\varphi_L$ to the central fiber: this is the main reason for introducing the additional singularities 
in the expression of the metric $\psi_{k, r}$.

We recall that we have $\displaystyle \varphi_{L_{r+1}}= \sum_{j\in J_{r+1}\subset J^\prime}\log |f_j|^2$, and hence the finiteness of the integral above is a consequence of the integrability lemma B.13, which will be stated and proved at the end of the present section. Thus, all the hypothesis required by the extension theorem A.0 are fulfilled, so there exists a section 
$U^{(kp+r+1)}_i$ extending $u^{(kp+r+1)}_i$, and which is $L^2$ with 
respect to $\varphi_E$. We remark that by induction we have 
$$\varphi^{(kp+ r)}\leq \sum_{j\in J^1}\log |f_{Y_j}|^2$$
and this is also the case for $\displaystyle{1\over p}\varphi_L$ (by definition), so we derive a similar conclusion for $\varphi_E$, given the expression (29). Thus, the case $r\leq p-2$ is
completely settled.

\smallskip

\noindent $\bullet$ We assume next that we have $r= p-1$. The section to be extended
during this step is say $u^{(kp+p)}_i\in V_{k+1, 0}$, 
so it verifies the following $L^2$ condition 
$$\int_{\cX_0}
|u^{(kp+p)}_i|^{2}e^{-\psi_{k+1,0}}< \infty \leqno (32)$$
hence we get
$$\int_{\cX_0}
{|u^{(kp+p)}_i|^{2}\over \prod _{j\in J\setminus J^\prime}|\sigma_j|^{2p\nu^j}}
e^{-\psi_{k, p-1}}d\lambda< \infty \leqno (33)$$
because of inequality (7); we denote by $\sigma_j$ the canonical section associated to the hypersurface $Y_j$. From the finiteness of the previous integral, we
derive two conclusions. In the first place, the section 
$$v^{(p)}:= {u^{(kp+p)}_i\over \prod _{j\in J\setminus J^\prime}\sigma_j^{p\nu^j}}$$
is {\sl holomorphic}. Secondly, $v^{(p)}$ belongs to the multiplier ideal sheaf 
$\displaystyle \cI(\psi_{k, p-1})$. By global generation property (22) of the bundle $A$, we therefore obtain
$${|u^{(kp+p)}_i|^{2}\over \prod _{j\in J\setminus J^\prime}|\sigma_j|^{2p\nu^j}}
\leq C\sum_j |u^{(kp+p-1)}_j|^2.\leqno(34)$$
Next, we write
$$(k+1)(pK_{\cX/\bD}+ L)+ A= K_{\cX/\bD}+ L_p+ kpK_{\cX/\bD}+ kL+ 
L^{(p-1)}+ A,$$
so we consider the bundle 
$$E:= L_p+ kpK_{\cX/\bD}+ kL+ 
L^{(p-1)}+ A.$$
In order to endow it with a metric, we recall that the Chern class of $L_p$ contains the current
$$p(\alpha+ \sum _{j\in J\setminus J^\prime}\nu^j[Y_j])+ \sum_{j\in J_p}
[Y_j].$$
The metric whose associated curvature form
is equal to the the first term of the previous sum is denoted by 
$\varphi_{p}^1$, and we define $\displaystyle \varphi_p:= \sum_{j\in J_p}\log |f_{Y_j}|^2$. They induce a metric on the bundle $E$ as follows
$$\varphi_E:= \varphi_{p}^1+ (1-\delta)\varphi_p+ 
\delta\wt\varphi_p+ (1-\varepsilon)\varphi^{(kp+ p-1)}+ 
\varepsilon\big((kp-1)\varphi_{\rm min, \tau}+ {1\over p}\varphi_L+ 
\wt\varphi_{p-1, A}\big).\leqno(35)$$
Its curvature current is positive as soon as $\delta\ll\varepsilon$ and $(kp- 1)\tau\ll 1$.

The $L^2$-norm of the section $u^{(kp+p)}_i$
with respect to $\varphi_E$ is {\sl finite}, provided that we have 
$$\int_{\cX_0}|u^{(kp+p-1)}_j|^{2\varepsilon}e^{-(1-\delta)\varphi_{L_{p}}-
\varepsilon(kp-1)
\varphi_{\rm min, \tau}-{\varepsilon\over p}\varphi_L}< \infty\leqno(36)$$
for each $j$.
As in the preceding case, the $L^2$ requirement reduces to (36) thanks to the 
inequality (34) above. The inequality (36) was already established during the analysis of the preceding case (31). 

In conclusion, there exists an extension $U^{(kp+p)}_i$
of the section $u^{(kp+p)}_i$, which moreover is integrable with respect to $\varphi_E$. This implies that the section $U^{(kp+p)}_i$ vanishes as required in 
lemma B.3, so the proof is complete. 
\hfill\qed

\medskip 
\noindent The next statement is a summary of the preceding considerations.

\claim B.4 Lemma|We have
$$\sup_j|u^{(kp)}_j(x)|^{{2\over kp}}\leq 
C^{-1}e^{{1\over kp}\varphi^{(kp)}_{\cX/\Delta}(x)}.
\leqno (37)$$
for any $k$, as well as for any $x\in \cX_0$. \hfill\qed
\endclaim

\smallskip
\proof. Indeed, we specialize the relation (24) for $\tau:= 0$, and $U:= U^{(kp)}_j$
(cf. lemma B.3); combined with lemma B.2, it gives the inequality (37) above. \hfill\qed

\medskip 
\noindent As a consequence of the regularization theorem due 
to J.-P. Demailly in [12], we have the 
following very precise estimate.
\claim B.5 Lemma|There exists a constant $C$ such that we have 
$${(kp-1)}\varphi_{0}(x)+ \sum_{j\in J}(1+\nu^j)
\log |f_{Y_j}|^2\leq C\log k+ \log\sup_j|u^{(kp)}_j(x)|^{{2}}$$
for any $x\in\cX_0$ and $k\in \bZ_+$ large enough. \hfill\qed
\endclaim 

\smallskip
\proof. We refer to the articles [12], [15]; the preceding inequality is obtained from the proof of the main theorem.   \hfill\qed

\medskip 
\proof \hskip 1mm (of 0.2). By theorem B.1, (b) we infer the 
existence of a positively curved limit metric
$$\varphi^{(\infty)}_{\cX/\Delta}:= \lim\sup_k {1\over kp}\varphi^{(kp)}_{\cX/\Delta}$$
on the $\bQ$-bundle $K_{\cX/\bD}+ 1/pL$.
By lemmas B.2-B.5, the metric $\varphi^{\infty}_{\cX/\Delta}$ is less singular than
$\varphi_0$ when restricted to the central fiber $\cX_0$. This metric is clearly more singular than $\varphi_{\rm min}$,
so the inequality (9) of theorem 0.2 is established (the uniformity of the constant $C$ in (9) is obtained by inspection of the proof of B.2-B.5). \hfill\qed

\claim B.6 Remark|{\rm The ``traditional" method of proving this kind of 
results {\sl does not seem to 
work} in this generalized setting. The reason is that we have to change the parameters $\delta, \varepsilon$ in the proof of $\rm B.3$ as $k\to \infty$, and the usual``concavity of the log" (in [41], [35]) cannot be applied in order to obtain the estimates needed for the justification of the limit metric above. The asymptotic $kp$-Bergman metric somehow converts the {\sl qualitative} information of the lemma B.3 into an {\sl effective estimate}.\hfill\qed}
\endclaim 

\noindent We state the second part of theorem 0.2 as a separate corollary.
\claim B.7 Corollary|Let $u$ be a section of the bundle $\displaystyle pK_{\cX_0}+ L_{|\cX_0}$, whose divisor of zeroes contains $p\sum_{j\in J^\prime}\rho^j_{\rm min, \infty}Y_{j0}+ p\sum_{j\in J\setminus J^\prime}\nu^jY_{j0}$; moreover, we assume that
$$\int_{\cX_0}|u|^{2\over p}e^{-{1\over p}\varphi_L}< \infty.\leqno(42)$$
Then there exists a section $U$ of $pK_{\cX}+ L$ extending $u$, and such that
$$\int_{\cX}|U|^{2\over p} e^{-{1\over p}\varphi_L}\leq C_0 \int_{\cX_0}
|u|^{2\over p} e^{-{1\over p}\varphi_L}.
$$
\hfill\qed
\endclaim
\noindent This statement can be seen as a generalization of [16], [22] where $L$ has an ample component. Certainly the convergence of the integral above
just means that $u$ vanishes on the log canonical part of 
$\displaystyle{1\over p}L$, but we prefer this formulation because it is very well adapted for the study of similar results under more general boundaries $L$. \hfill\qed

\proof. The vanishing properties of $u$ together with the inequality (10) shows the existence of some constant $C$ such that we have
$$|u|^2e^{-p\varphi_{\rm min}}\leq C< \infty$$
on the central fiber. We have $\displaystyle pK_{\cX}+L= K_{\cX}+ {p-1\over p}(pK_{\cX}+L)+ {1\over p}L$, and the above inequality shows that the $L^2$ norm of $u$ with respect to the metric $(p-1)\varphi_{\rm min}+ {1\over p}\varphi_L$ is finite (we use here the $L^{2\over p}$
convergence hypothesis (42) in B.7). The proof ends thanks to the Ohsawa-Takegoshi theorem. \hfill\qed

\medskip 

\noindent In order to prove the inequality (10) stated in the introduction, we will assume that
$$\nu^j< 1\leqno(38)$$
for all $j\in J$. 

Let $\displaystyle L^\prime:=  L_{|\cX_0} - p\sum_{j\in J^\prime}\rho^j_{\rm min, \infty}Y_{j0}-
p\sum_{j\in J\setminus J^\prime}\nu^jY_{j0}$; by definition of the set $J^\prime$ we see that it is a pseudoeffective 
$\bR$-divisor on
the central fiber $\cX_0$, whose adjoint 
$$pK_{\cX_0}+ L^\prime\leqno(39)$$
is pseudoeffective. This property is a consequence of the fact that $\displaystyle \varphi_{\rm min|\cX_0}$
is well defined, so that we can write 
$$\Theta_{\rm min|\cX_0}= \sum_{j\in J}\rho^j_{\rm min}[Y_{j0}]+ \Lambda_0;$$
the observation is that $\rho^j_{\rm min}\geq \rho^j_{\rm min, \varepsilon}$ for any 
$\varepsilon> 0$, and thus the same inequality holds for the limit.

We denote by $\varphi_{\rm min}^\prime$
a metric with minimal singularities corresponding to the bundle (39); a direct consequence of the theorem 0.2 is the next statement.

\claim B.8 Corollary|We have
$$\big|\varphi_{\rm min|\cX_0}- \sum_{j\in J^\prime}\rho^j_{\rm min \infty}\log|f_j|^2- 
\sum_{j\in J\setminus J^\prime}\nu^j\log|f_j|^2 - \varphi_{\rm min}^\prime\big|\leq C$$
pointwise on $\cX_0$. \hfill\qed
\endclaim
\proof. We first observe that the expression
$$\psi:= \varphi_{\rm min|\cX_0}- \sum_{j\in J^\prime}\rho^j_{\rm min, \infty}\log|f_j|^2- 
\sum_{j\in J\setminus J^\prime}\nu^j\log|f_j|^2$$ corresponds to a 
positively curved metric of the bundle in (39) (despite of the minus signs in its definition...), and thus we have
$$\psi\leq  \varphi_{\rm min}^\prime+ C$$
by definition of the minimal metric associated to a cohomology class.

In the opposite sense, we note that we have
$$\varphi^{\infty}_{\cX/\bD|\cX_0}\geq \sum_{j\in J^\prime}\rho^j_{\rm min}\log|f_j|^2+ 
\sum_{j\in J\setminus J^\prime}\nu^j\log|f_j|^2+ \varphi_{\rm min}^\prime + C\leqno (40)$$
by inequality (10), where the metric $\varphi_0$ corresponds to the right hand side of the above relation. Finally, we clearly have
$$\varphi_{\rm min}\geq \varphi^{\infty}_{\cX/\bD}\leqno (41)$$
and the corollary $\rm B.7$ is proved. \hfill\qed
\smallskip

\medskip
\noindent We turn now to the analysis of the {\sl compact version of 0.2}. 
This means that we replace the couple $(\cX, \cX_0)$ by $(X, S)$, where $X$ be a projective manifold, and $S\subset X$ is a non-singular hypersurface. We consider a line bundle $L\to X$ such that 
$$p([\Delta]+ \alpha)\in c_1(L)\leqno(43)$$
where $\Delta= \sum_{j\in J}\nu^jY_j$ and $\alpha\geq 0$ have the properties $\bf (a), \bf (b)$ and $(\dagger)$ in the introduction.

The class $\{K_X+ S+ 1/pL\}$ is assumed to be pseudoeffective, and we denote by $\varphi_{\rm min}$ a metric with minimal singularities corresponding to it; the associated curvature current will be denoted by
$\Theta_{\rm min}$. 
As before, we suppose that 
$$\nu_{\rm min}(\{K_X+ S+ 1/pL\}, S)= 0\leqno(44)$$ 
and then we 
can define the quantities $\rho^j_{\rm min, \infty}$ exactly as in 
the previous case: let $\Theta_{\rm min, \varepsilon}$ be a current with minimal singularities within the class $\{K_X+ S+ 1/pL+ \varepsilon A\}$; we have
$$\Theta_{\rm min, \varepsilon|S}:= \sum_{j\in J} \rho^j_{\rm min, \varepsilon}[Y_{j|S}]+ \Lambda_{S, \varepsilon}\leqno(45)$$
(thanks to the assumption (44) above) where $(\rho^j_{\rm min, \varepsilon})$ are positive real numbers, and $\Lambda_{S, \varepsilon}$
is a closed positive current defined on $S$. The limit of  $\rho^j_{\rm min, \varepsilon}$ is denoted by $\rho^j_{\rm min, \infty}$. With this quantities we define the set $J^\prime\subset J$ as in the introduction, and we assume that $\bf (c)$
holds as well.

Another part of the data is a {\sl positively curved metric} 
$\varphi_0$ on the bundle $K_X+ S+ 1/pL_{|S}$; we assume that it satisfies the properties (7) and (8) in the introduction.
We discuss next 
the following version of the theorem 0.2.

\claim B.9 Theorem|Under the hypothesis above, the metric $\varphi_{\rm min|S}$ is not identically $-\infty$ and we have 
$$\varphi_{\rm min|S}\geq C+ \varphi_0$$
pointwise on $S$. \hfill\qed
\endclaim

\proof. The following arguments are completely similar to the ones provided for the proof of 0.2 along the steps B.2-B.5. We will explain next the few things which are to be changed in order to conclude.

\smallskip

\noindent ${\bullet}$ We consider the space  
$$V_{k, r}:= H^0\Big(S,  \big(kp(K_{X}+ S)+ kL+ L^{(r)}+ A_{|S}\big)\otimes \cI(\psi_{k, r})\Big)\leqno(46)$$
where $\displaystyle \psi_{k, r}:= (kp-1)\varphi_{0}+ \sum_{j\in J}(1+ \nu^j)\log|f_{Y_j}|^2+ \wt \varphi_{r, A}$; 
is a positively curved metric.
Let $(u^{(kp+ r)}_j)$ be an orthonormal basis of the space $V_{k, r}$; then we have
$$\int_{S}
|u^{(kp+r)}_j|^{2/kp}\exp \big(-{{1}\over {p}}\varphi_{L}-{{1}\over {kp}}\varphi_{r, A}\big)\leq C\leqno (47)$$
for all $k\geq 1$ thanks to the H\"older inequality, see B.2 (in the expression under the integral sign, we identify 
$u^{(kp+r)}_j$ with a section of the bundle 
$kpK_{S}+ kL+ L^{(r)}+ A_{|S}$). Given the singularities of the metric $\displaystyle \psi_{k, r}$, we infer that the section $u^{(kp+r)}_j$
vanishes along the divisor $\sum_{m\in J}Y_{m|S}$.

\smallskip

\noindent ${\bullet}$
The algorithm used in the proof of the lemma B.3 shows that given the integers $(k, r, j)$, the 
corresponding section $u^{(kp+ r)}_j$ admits {\sl some} extension to $X$. In fact, we show next that we can construct an extension of $u^{(kp+ r)}_j$
which verifies an {\sl effective estimate}, crucial for the rest of the proof.

Let $J^1\subset J$ be the set of indexes $j\in J$ such that $\nu^j= 1$. Exactly as in the proof of B.3 we show that {\sl there exists an extension 
$\wt U^{(kp+ r)}_j$ of the section $u^{(kp+ r)}_j$ which vanishes on the divisor
$\sum_{j\in J^1}Y_j$}. In particular, the section $\wt U^{(kp+ r)}_j$ verifies the inequality
$$\int_{X}
|\wt U^{(kp+r)}_j|^{2\over kp}\exp \big(-{k\varphi_{L}+ \varphi_{r, A}\over kp}-\varphi_S\big)< \infty\leqno(48)$$
The next claim is that we can choose an extension $U^{(kp+ r)}_j$
of the section $u^{(kp+ r)}_j$ such that 
$$\int_{X}
|U^{(kp+r)}_j|^{2\over kp}\exp \big(-{k\varphi_{L}+ \varphi_{r, A}\over kp}-\varphi_S\big)\leq C_0\int_{S}|u^{(kp+r)}_j|^{2\over kp}
\exp \big(-{k\varphi_{L}+ \varphi_{r, A}\over kp}\big)\leqno (49)$$
for all $k\gg 0$, where the constant $C_0$ is independent of $k$. Indeed, the proof of the $L^{p}$-Ohsawa-Takegoshi in [4] shows that the extension of $u^{(kp+ r)}_j$ which minimize the left hand side of (49) will automatically verify the estimates; we do not provide here any further details, but rather remark that this is the only place in the proof where the condition $(\dagger)$ is used.



\smallskip
\noindent ${\bullet}$ We denote by $\psi^{(kp)}$ the metric on the bundle $kp(K_X+ S)+ kL+ A$ associated to the set of sections $U^{(kp)}_j$. The inequality (49) above shows that 
$$\psi^{(\infty)}:= \lim_{\rm reg}\sup_{k\to \infty}{1\over kp}\psi^{(kp)}\leqno(50)$$
is a positively curved metric on the $\bQ$-bundle $K_X+S+ {1\over p}L$. 
Its restriction to $S$ is greater than $\varphi_0$ (up to a constant), by the same arguments as 
in B.4 and B.5. Certainly the metric $\psi^{(\infty)}$ is smaller than the metric with minimal singularities, so that the theorem 
B.8 is proved.
\hfill\qed
\medskip

\claim Question|{\rm One of the important points in the proof of theorem B.8 was that we can construct the extensions $U^{(kp)}_j$
vanishing along $\sum_{j\in J_1}Y_j$. We describe next a related question.

Let $\Delta$ be an effective divisor on $X$; one can define
a notion of $\Delta$-minimal metric on $K_X+S+ {1\over p}L$, i.e.
the upper envelope of all normalized, positively curved metrics which are {\sl at least} as singular as the quasi-psh function associated to the divisor 
$\Delta$. We denote this object by $\varphi_{\rm min, \Delta}$ and we assume that its restriction to $S$ is not identically $-\infty$. The question is to identify the restriction $\varphi_{\rm min, \Delta|S}$; more precisely, we ask for a criteria similar to the corollary B.7. Unfortunately, the methods used and developed in this article do not seem to be very helpful in this direction. \hfill\qed}
\endclaim

\medskip
\claim B.10 Remark|{\rm A large part of the proof of B.9 can be applied in a 
more general setting, but it only gives a qualitative result.

For example, instead of the set $J^\prime$ in the introduction we can define
$$J_0^\prime:= \{j\in J : \rho^j_{\rm min, \infty}< \nu^j\leq 1\};$$
we also formulate the next condition.

\noindent ${\bf (c^\prime)}$ 
For any subset $I\subset J^\prime_0$ and for any $\varepsilon> 0$ 
the restriction of the current 
$\Lambda_{0 \varepsilon}$ defined in the relation (5) to the intersection $\displaystyle \bigcap_{m\in I}Y_m\cap S$
is well-defined.\hfill\qed

\noindent The hypothesis are the same as for B.9, except that we use $J_0^\prime$ 
and $\bf (c^\prime)$ to replace $J^\prime$ and respectively $\bf (c)$.
Then we infer the following result, which will be crucial for the next section.

Consider the space $V_{k, r}$ as in (46); then relation (47) still holds without further modifications. Concerning the second bullet in the proof of B.8, we can only show the existence of an extension $U^{(kp+r)}_j$ of $u^{(kp+r)}_j$, without the estimate (48). The family $\displaystyle (U^{(kp)}_j)_{j}$
defines a metric $\psi^{(kp)}$ on the bundle $kp(K_X+ S)+ kL+ A$, and we have the estimate
$$\psi^{(kp)}\geq {(kp-1)}\varphi_{0}+ \sum_{j\in J}(1+\nu^j)
\log |f_{Y_j}|^2+ C$$
pointwise on $S$. A last observation is that if the bundle $p(K_X+S)+ L$
happens to be big, then we can define the space $V_{k, r}$ by using sections of $kp(K_X+S)+ kL+ L^{(r)}_{|S}$ (i.e. without the additional twisting with $A$) by a slight modification of the weight as in [15], section 17. Then the family of extensions $\displaystyle (U^{(kp)}_j)_{j}$ are sections of 
$kp(K_X+S)+ kL+ L^{(r)}$, and the above inequality becomes
$$\psi^{(kp)}\geq {\big((k-k_0)p-1\big)}\varphi_{0}+ \psi+ \sum_{j\in J}(1+\nu^j)
\log |f_{Y_j}|^2+ C,\leqno(\bu)$$
where $\psi$ is a metric on on $k_0p(K_X+S)+ k_0L$ whose curvature current is greater than a K\"ahler metric, and $k_0$ is a large enough integer.}

\endclaim

\medskip
A consequence of the corollary B.9 is that the norm with respect to the metric $\displaystyle e^{-\varphi_{\rm min|S}}$
of any section $u$ of the bundle
$p(K_X+S)+ L_{|S}$ whose zero set contains the divisor $\sum_{j\in J^\prime}\rho^j_{\rm min, \infty}Y_{j|S}+ \sum_{j\in J\setminus J^\prime}\nu^jY_{j|S}$
is {\sl uniformly bounded}. We show in the next corollary that $u$ admits an extension to $X$ in the klt case. 

\claim B.11 Corollary|In addition to the hypothesis in B.9, we assume that $\nu^j< 1$ for all $j\in J$. Then any section 
$$u\in H^0\big(S, p(K_X+S)+ L_{|S}\big)$$
whose zero set contains the divisor $\sum_{j\in J^\prime}\rho^j_{\rm min}Y_{j|S}+ \sum_{j\in J\setminus J^\prime}\nu^jY_{j|S}$ admits an extension to $X$. \hfill\qed
\endclaim

\proof. Let $A\to X$ be an ample enough line bundle, such that 
$$u\otimes \sigma_A$$
extends to $X$, where $\sigma_A$ is a non-zero section of $A$. We denote by $U_A$ the corresponding extension, and we consider the section
$$u^{\otimes 2}\otimes \sigma_A \in H^0\big(S, 2p(K_X+S)+ 2L+ A_{|S}\big)$$
We will 
construct next an of extension of $u^{\otimes 2}\otimes \sigma_A$ which 
is divisible by $U_A$; the quotient will be the desired extension of $u$.

To this end, we have the equality
$$2p(K_X+S)+ 2L+ A= K_X+ S+ p(K_X+S)+ L+ A+ (p-1)(K_X+S+ {1\over p}L)+ {1\over p}L,$$
and for each positive $\varepsilon$, we consider the metric
$$ \log |U_A|^2+ (p-1)\varphi_{\rm min}+ {1\over p}\varphi_L.$$
Its curvature form is both semi-positive, and it dominates 
$\displaystyle {\gamma\over p} \Theta _h\big(\cO(S)\big)$. In order to apply [13], we have to check next the integrability condition  
$$\int_S{|u^2\otimes\sigma_A|^2 \over |u\otimes\sigma_A|^2}
\exp\big(-(p-1)\varphi_{\rm min}- {1\over p}\varphi_L\big)<\infty.$$
This is however obvious, since by corollary $\rm B.8$ we have
$$|u|^2\leq C e^{p\varphi_{\rm min|S}}$$
at each point of $S$, and the restriction to $S$ of $\displaystyle e^{-1/p\varphi_L}$ is convergent.

Thus, we obtain a section 
$$V \in H^0\big(X, 2p(K_X+S)+ 2L+ A\big)$$
whose restriction to $S$ is equal to $u^2\otimes\sigma_A$, and such that
$$\int_X{|V|^2 \over |U_A|^2}
\exp\big(-(p-1)\varphi_{\rm min}- {1\over p}\varphi_L\big)<\infty$$
This in turn implies that the quotient $\displaystyle {V\over U_A}$ is a holomorphic section of $p(K_X+S)+ L$, and it is equal to $u$ when restricted to $S$. \hfill\qed

\medskip
\noindent We establish next the following 
integrability criteria, which was used several times in this paragraph.

\claim B.13 Lemma|Let $\Theta$ be a closed (1,1)--current on a K\"ahler manifold
$(X, \omega)$, such that 
$$\Theta\geq -C\omega$$
for some positive constant $C$. Let $\displaystyle (Y_j)_{j= 1,..., r}$ be a finite set of hypersurfaces, which are supposed to be non-singular and to have normal crossings. Moreover, we assume that 
the restriction of $\Theta$ to the intersection $\cap _{i\in I}Y_i$ is well defined,
for any $I\subset \{1,..., r\}$. Then there exists a positive $\varepsilon_0= \varepsilon_0(\{\Theta\}, C)$ depending only on the cohomology class of the current and on the 
lower bound $C$ such that for any $\delta\in ]0, 1]$ and $\varepsilon \leq \varepsilon_0$ we have

$$\int_{(X, x)}\exp \big(-(1-\delta)\sum_{j=1}^r\log|f_{Y_j}|^2-\varepsilon\varphi_\Theta\big)d\Lambda< \infty,$$
at each point $x\in X$. \hfill\qed
\endclaim

\proof. Let $x\in X$ be an arbitrary point; 
we can assume that $x\in Y_1\cap...\cap Y_b$ and $y\not\in Y_k$ for some $b\leq r, k\geq b+1$.
For each $p=1,..., b$ we define the complete intersection
$$\Xi_p:=  Y_1\cap...\cap Y_p$$
and the Skoda integrability theorem imply that 
$$\int_{(\Xi_b, x)}\exp (-\varepsilon\varphi_\Theta)d\lambda< \infty$$
for any $\varepsilon\leq \varepsilon_0\ll 1$ (we remark that here we use the hypothesis 
concerning the restriction of $\Theta$ to the sets $\Xi_p$). 

By the local version of the Ohsawa-Takegoshi theorem (cf. [13]), we can extend the constant function equal to 1 on $\Xi_b$ to a holomorphic function $f_{b-1}\in \cO(\Xi_{b-1}, x)$ such that
$$\int_{(\Xi_{b-1}, y)}|f_{b-1}|^2\exp \big(-(1-\delta)\log|f_{Y_{b}}|^2- \varepsilon \varphi_\Theta\big)d\lambda< \infty;$$
we repeat this procedure $b$ times, until we get a function $f_0\in \cO(X, x)$
such that
$$\int_{(X, x)}|f_0|^2\exp \big(-(1-\delta)\sum_{j=1}^r\log|f_{Y_j}|^2-\varepsilon\varphi_\Theta\big)d\lambda< \infty.$$
Since the function $f_0$ is constant equal to 1 in a open set centered at $y$ in $\Xi_b$,
we are done, except for the uniformity of $\varepsilon_0$. 

Indeed, the fact that the $\varepsilon_0$ only depends on the quantities in the above statement is a consequence of the fact that the Lelong numbers of closed positive currents on K\"ahler manifolds are bounded by the cohomology class of the current.\hfill\qed

As a side remark, one can see that the preceding statement hold true under the 
weaker assumption 
$$\nu_{\cap _{i\in I}Y_i}(\Theta)= 0$$
that is to say, we claim that the previous lemma is true if the generic Lelong number of $\Theta$ along all the intersections above is zero. Indeed, one can apply the regularization theorem 
stated in 2.1 combined with the H\"older inequality in order to derive the general result; we leave the details to the interested reader.
\hfill\qed


\bigskip

\section {\S C. Further applications}

We will prove in this paragraph an extension statement which was used in the 
article [4] (theorem B.1.2). We first recall the general set-up in 
[4] (and use the notations in that article).

Let $X$ be a {\sl normal projective variety}, and let $\Delta$ be an effective Weil $\bQ$-divisor
on $X$, such that $K_X+ \Delta$ is $\bQ$-Cartier. 
We consider $W\subset X$ an {\sl exceptional center} of $(X, \Delta)$; in other words, we assume that
there exists a log-resolution $\mu: X^\prime\to X$ of the pair $(X, \Delta)$ together with a decomposition of the inverse image of the $\bQ$-divisor $K_X+ \Delta$ as follows
$$\mu^\star(K_X+ \Delta)=  K_{X^\prime}+ S+ \Delta^\prime+ R- \Xi,\leqno (\dagger)$$
such that:
{
\itemindent 3.5mm
\smallskip
\item {$\bullet$} $S$ is an irreducible hypersurface, such that $W= \mu(S)$;
\smallskip
\item {$\bullet$} $\Delta^\prime:= \sum_j a^jY_j$, where $x_0\in \mu(Y_j)$ and $a^j\in ]0, 1[$;
\smallskip
\item {$\bullet$} The divisor $R$ is effective, and a hypersurface $Y_j$ belongs to 
its support if either $x_0\not\in \mu(Y_j)$, or $Y_j\cap S= \emptyset$ (so in particular
the restriction $R_{|S}$ is $\mu_{|S}$-vertical);

\smallskip
\item {$\bullet$} The divisor $\Xi$ is effective and $\mu$-contractible; in addition, we assume that the support of the divisors of the right hand side of the formula $(\dagger)$ 
has strictly normal crossings.

}

\noindent In general, the center $W$ is singular, and we will assume that the restriction map
$\displaystyle \mu_{|S}: S\to W$
factors thru the desingularization $g: W^\prime\to W$, so that we have
$$\mu_{|S}= g\circ p$$
where $p: S\to W^\prime$ is a surjective projective map. 

Before stating our next result, we introduce a last piece of notation: let $A$ be an ample bundle on $X$, and let $F_1,..., F_k$ be a set of smooth hypersurfaces of $W^\prime$ with strictly normal crossings, such that there exists positives rational numbers 
$(\delta_j)$ for which the $\bQ$-bundle
$$g^\star(A)-\sum \delta^jF_j+ \varepsilon K_{W^\prime}\leqno (51)$$
is semi-positive (in metric sense) for any $\varepsilon$ small enough, 
and such that $g(F_j)\subset W_{\rm sing}$ for each $j$.
Indeed a set $(F_j)$ with the properties specified above does exists, see e.g. [4].

The family $(F_j)$ induces a decomposition of the divisor $\Xi$ as follows
$$\Xi= \Xi_1+ \Xi_2$$
where by definition $\Xi_{1}$ is the part of the divisor $\Xi$
whose support restricted to $S$ is  
mapped 
by $p$ into $\cup_j F_j$.
\smallskip

\noindent The result we will prove next is the following.
\claim C.1 Theorem|Let $T$ be any  closed positive $(1,1)$-current, such that there exists a line bundle $E$ on 
$X$ with the property that
$T\equiv \mu^\star (E)_{|S}+ m(K_{S/{W^\prime}}+ \Delta^\prime_{|S})$. Then we have 
$$T\geq m[\Xi_{2|S}]$$ 
in the sense of currents on $S$.
\hfill\qed

\endclaim
 
 \proof. We start with a few reductions; by hypothesis (51) we infer the existence of a non-singular and semi-positive $(1,1)$-form $\alpha$ such that 
 $$\alpha\equiv C\big(g^\star (A)- \sum_i\delta^iF_i\big)+ mK_{W^\prime}$$
 where $C> 0$ is a large enough 
 constant. Therefore we obtain
 $$T+  C\sum_i\delta^i[p^\star(F_i)]+ p^\star(\alpha)\equiv 
 {\mu^\star\big(E+ CA\big)+ m(K_{X^\prime}+ S+ \Delta^\prime)}_{|S}.
 \leqno(52)$$
By definition of the decomposition $\Xi= \Xi_1+ \Xi_2$, irreducible components of the 
divisors $\sum_i\delta^ip^\star(F_i)$ and respectively $\Xi_{2|S}$ are {\sl disjoint}
 (we use here the fact that $\mu$ is a log resolution). Therefore, in order to prove the previous theorem it is enough to show that the current
$$\Theta:= T+  C\sum_i\delta^i[p^\star(F_i)]+ p^\star(\alpha)\leqno (53)$$
verifies the inequality
$$\Theta\geq m[\Xi_{2|S}]\leqno(54)$$ 
in the sense of currents on $S$.

\noindent We show next that the inequality (54) is a consequence of remark B.10.
By relation $(\dagger)$, we have
$${\mu^\star\big(E+ CA+ m(K_X+ \Delta)\big)+ m\Xi}= {\mu^\star\big(E+ CA\big)+ m(K_{X^\prime}+ S+ \Delta^\prime+ R)}.
 \leqno(55)$$
 We can assume that $A$ is large enough, such that 
 $\displaystyle K_X+ \Delta+ {1\over m}(CA+ E)$ is ample as well; then we claim that {\sl the metric with minimal singularities $\varphi_{\rm min}$ associated to the big class 
$${1\over m}\{\mu^\star\big(E+ CA+ m(K_X+ \Delta)\big)+ m\Xi\}.
 \leqno(56)$$
 has the same singularities as $\varphi_{\Xi}$}. 
 
 To see this, we first recall that by Hartogs principle the zeroes divisor of any section of the bundle
 $$k\big(\mu^\star\big(E+ A+ m(K_X+ \Delta)\big)+ m\Xi\big)$$
 is greater than $km\Xi$
 (indeed, the Hartogs principle is still valid even if $X$ may have singularities; the thing is that it is normal). Since the cohomology class (56) is big, the regularization of closed positive currents 
 in [12] shows that  
{the} positively curved metric $e^{-\varphi_{\rm min}}$ can be approximated with holomorphic sections of the bundle
$$k\big(\mu^\star\big(E+ A+ m(K_X+ \Delta)\big)+ \Xi\big)$$
in a very precise way: there exists a (singular) metric $\psi$ on the bundle
$k_0\big(\mu^\star\big(E+ A+ m(K_X+ \Delta)\big)+ \Xi\big)$ such that 
$(k-k_0)\varphi_{\rm min}+ \psi$ is smaller than the metric induced by a family of holomorphic sections of the bundle above, see [12].
As $k\to\infty$, we have proved our claim.

Let $\varphi_{\Theta+ R}$ be the metric of the bundle 
$$K_{X^\prime}+ S+ \Delta^\prime+ R+ {1\over m}\mu^\star(E+ A)_{|S}$$
whose curvature current is equal to $\Theta+ [R]$.

We apply the results obtained in B.10 for $\varphi_0:= \varphi_{\Theta+ R}$.
The hypothesis needed for the inequality $(\bu)$ to hold are satisfied, given the structure of $\varphi_{\rm min}$ discussed above. Hence, despite the fact that in the inequality $(\bu)$ we do not have any upper bound for the algebraic metric
$\psi^{(kp)}$, we can infer that the the metric $\varphi_{\Theta+ R}$
is more singular than $\displaystyle m\varphi_{\Xi_2}$.
The proof of theorem C.1 is finished.\hfill\qed


\vskip 15pt

\vfill \eject 
\section{References}
\bigskip

{\eightpoint

\bibitem [1]&Berndtsson, B.:& On the Ohsawa-Takegoshi extension theorem;& Ann.\ Inst.\ Fourier (1996)&

\bibitem [2]&Berndtsson, B.:& Integral formulas and the Ohsawa-Takegoshi extension theorem;& Science in Chi\-na, Ser A Mathematics,  2005, Vol 48&

\bibitem [3]&Berndtsson, B., P\u aun, M.:& Bergman kernels and the pseudo-effectivity of the relative canonical bundles;& arXiv:math/0703344&

\bibitem [4]&Berndtsson, B., P\u aun, M.:& Bergman kernels and subadjunction type results,&\ revised version of arXiv:08043884&

\bibitem [5]&Birkar, C., Cascini, P., Hacon, C., McKernan, J.:&\ Existence of minimal models for varieties of log general type~;& to appear in J.\ Amer.\ Math.\ Soc&

\bibitem [6]&Boucksom, S.& C\^ones positifs des vari\'et\'es complexes compactes,& Thesis, Grenoble 2002&

\bibitem [7]&Boucksom, S.& Divisorial Zariski decompositions on compact complex manifolds;&  Ann. Sci. Ecole Norm. Sup. (4)  37  (2004),  no. 1, 45--76&

\bibitem [8]&Campana, F.:& Special Varieties and classification Theory;&
Annales de l'Institut Fourier 54, 3 (2004), 499-665&
 
\bibitem [9]&Campana, F., Peternell, Th:& Geometric stability of the cotangent bundle and the universal cover of a projective manifold;& arXiv:math/0405093& 

\bibitem [10]&Claudon, B.:& Invariance for multiples of the twisted canonical bundle;& math.AG/0511736, to appear in Annales de l'Institut Fourier&

\bibitem [11]&Demailly, J.-P.:&\ Singular hermitian metrics on positive line bundles,& Proc. Conf. Complex
algebraic varieties (Bayreuth, April 2Ð6, 1990), edited by K. Hulek, T. Peternell,
M. Schneider, F. Schreyer, Lecture Notes in Math., Vol. 1507, Springer-Verlag, Berlin,
1992&

\bibitem [12]&Demailly, J.-P.:&\ Regularization of closed positive currents and Intersection Theory; &J. Alg. Geom. 1 (1992) 361-409&

\bibitem [13]&Demailly, J.-P.:&  On the Ohsawa-Takegoshi-Manivel  
extension theorem;& Proceedings of the Conference in honour of the 85th birthday of Pierre Lelong, 
Paris, September 1997, Progress in Mathematics, Birkauser, 1999&

\bibitem [14]&Demailly, J.-P.:& K\"ahler manifolds and transcendental techniques in algebraic geometry;&  Plenary talk and Proceedings of the Internat. Congress of Math., Madrid (2006), 34p, volume I&

\bibitem [15]&Demailly, J.-P.:&\ Analytic methods in algebraic geometry,&
on the web page of the author, December 2009&

\bibitem [16]&Ein, L., Popa, M.:& Extension of sections via adjoint ideals;& private communication, June 2007, arXiv:0811.4290&

\bibitem [17]&de Fernex, T., Hacon, C.D.:&\ Rigidity properties of Fano varieties,& arXiv:0911.0504&
 
\bibitem [18]&de Fernex, T., Hacon, C.D.:&\ Deformations of canonical pairs and Fano varieties,&  
To appear in J. Reine Angew.\ Math,\ arXiv:0901.0389& 

\bibitem [19]&Hacon, C. D.:&\ Extension theorems and the existence of flips,& Series of 8 lectures given at
Oberwolfach Math. Institute, Oct. 12-18, 2008&

\bibitem [20]&Hacon, C.D., McKernan, J.:& Boundedness of pluricanonical maps of varieties of general type;&
 Invent.\ Math.\ Volume {\bf 166}, Number 1 / October, 2006, 1-25&
 
\bibitem [21]&Hacon, C.D., McKernan, J.,:& Extension theorems and the existence of flips,& In ``Flips for 3-
folds and 4-foldsÓ (ed. A. Corti), Oxford lecture series in mathematics and its applications
35, Oxford University Press, 2007&

\bibitem [22]&Hacon, C.D., McKernan, J.,:& Existence of minimal models for varieties of general type,
II, (Existence of pl-flips),& arXiv:math.AG/0808, 2009&

\bibitem [23]&Kawamata, Y.:& Pluricanonical systems on minimal algebraic varieties&\ Invent. Math.  79  (1985),  no. 3&

\bibitem [24]&Kawamata, Y.:& Subadjunction of log canonical divisors;& Amer.\ J.\ Math.\  {\bf 120} (1998) 893--899&

\bibitem [25]&Kawamata, Y.:& On the extension problem of pluricanonical forms,& Preprint Univ. of
Tokyo at Komaba, September 1998, Contemporary Math. 241(1999), 193Ð207&

\bibitem [26]&Kawamata, Y.:& Deformation of canonical singularities,& 
J.\ Amer. \ Math.\ Soc.,\ 12 (1999),
85Ð92&

\bibitem [27]&Kim, Dano.:& $L^2$ extension of adjoint line bundle sections;& arXiv:0802.3189, to appear in Ann.\ Inst.\ Fourier&

\bibitem [28]&Kim, Dano.:& Personal communication;& October 2009&

\bibitem [29]&Lazarsfeld, R.:& Positivity in Algebraic Geometry;& Springer, Ergebnisse der Mathematik und ihrer Grenzgebiete&
 
\bibitem [30]&Levine, M.:& Pluri-canonical divisors on K\"ahler manifolds;&  Invent. Math.  74  (1983),  no. 2, 293--303&

\bibitem [31]&Nakayama, N.:&Zariski decomposition and abundance~;&  MSJ Memoirs  {\bf 14}, Tokyo (2004)&
 
 \bibitem [32]&Ohsawa, T., Takegoshi, K.\ :& On the extension of $L^2$
holomorphic functions;& Math.\ Z.,
{\bf 195} (1987), 197--204&
 
 \bibitem [33]&Ohsawa, T.\ :& On the extension of $L\sp 2$ holomorphic functions. VI. A limiting case;&
Contemp.\ Math. \ ,  {\bf 332} (2003)
Amer. Math. Soc., Providence, RI&

 \bibitem [34]&Ohsawa, T\ :&Generalization of a precise L2 division theorem,& Complex analysis in several
variables Ð Memorial Conference of Kiyoshi OkaÕs Centennial Birthday, 249Ð261, Adv.
Stud. Pure Math., 42, Math. Soc. Japan, Tokyo, 2004&

\bibitem [35]&P\u aun, M.:& Siu's Invariance of Plurigenera: a One-Tower Proof;&preprint 2005, to appear in J.\ Diff.\ Geom&

\bibitem [36]&P\u aun, M.:& Relative critical exponents, non-vanishing and metrics with minimal singularities& arXiv:08073109&

\bibitem [37]&Shokurov, V.:& A non-vanishing theorem~;&\ Izv. Akad. Nauk SSSR
(49) 1985&

\bibitem [38]&Shokurov, V. V.:&\ Prelimiting flips& Tr. Mat. Inst.
Steklova, 240 (Biratsion. Geom. Linein. Sist. Konechno Porozhdennye
Algebry), 82Ð219, 2003&

\bibitem [39]&Siu, Y.-T.:& Analyticity of sets associated to Lelong numbers and the extension of closed positive currents& Invent. Math.  27  (1974), 53--156&
 
\bibitem [40]&Siu, Y.-T.:& Invariance of Plurigenera;& Inv.\ Math.,
{\bf 134} (1998), 661-673&

\bibitem [41]&Siu, Y.-T.:& Extension of twisted pluricanonical sections with plurisubharmonic weight and invariance of semipositively twisted plurigenera for manifolds not necessarily of general type;& Complex geometry (G\"ottingen, {\bf 2000}),  223--277, Springer, Berlin, 2002&

\bibitem [42]&Siu, Y.-T.:& Multiplier ideal sheaves in complex and algebraic geometry;& Sci.\ China Ser.  {\bf A 48}, 2005&

\bibitem [43]&Siu, Y.-T.:& A General Non-Vanishing Theorem and an Analytic Proof of the Finite Generation of the Canonical Ring;& arXiv:math/0610740&

\bibitem [44]&Takayama, S:& Pluricanonical systems on algebraic varieties of general type;&Invent.\ Math.\
Volume {\bf 165}, Number 3 / September, 2005, 551-587&

\bibitem [45]&Takayama, S:& On the Invariance and Lower Semi--Continuity
of Plurigenera of Algebraic Varieties;& J. Algebraic Geom.  {\bf 16 } (2007), no. 1, 1--18&

\bibitem [46]&Tsuji, H.:& Extension of log pluricanonical forms from subvarieties;& math.CV/0511342&

\bibitem[47]&Tsuji, H.:& Canonical singular hermitian metrics on relative canonical bundles;& prepublication
Sophia University, Japan, arXiv:math.AG/0704.0566&

\bibitem[48]&Tsuji, H.:&Canonical volume forms on compact K¬ahler manifolds,& prepublication Sophia
University, Japan, arXiv: math.AG/0707.0111&
 
\bibitem [49]&Varolin, D.:&  A Takayama-type extension theorem;&  math.CV/0607323, to appear in Comp.\ Math&
     
\bibitem [50]&Viehweg, E.:& Quasi-Projective
Moduli for
Polarized Manifolds;& Springer-Verlag, Berlin, Heidelberg, New York, 1995
as: Ergebnisse der Mathematik und ihrer Grenzgebiete, 3. Folge, Band 30&

}

\bigskip
\noindent
(version of February 20, 2010, printed on \today)
\bigskip\bigskip
{\parindent=0cm
Bo Berndtsson,  
bob@math.chalmers.se\\
Mihai P\u aun,
paun@iecn.u-nancy.fr
}

\end

\medskip
We recall the next form of the Hartogs extension principle, concerning the line bundles 
$$\wh L:= \mu^\star L+ \wh E$$
where $L\to X$ is a line bundle and $\wh E= \sum b^{(q)}E_q$ is effective and contractible.
The map $\mu: \wh X\to X$ is assumed to be a composition of blow-up maps with non-singular centers. 

\claim Hartogs Lemma|Let $\wh T\in c_1(\wh L)$ be a closed positive current. Then 
$$\wh T\geq [\wh E];$$
in particular, the zero set of any holomorphic section of some line bundle numerically equivalent to 
$\wh L$ contains the divisor $\wh E$.
\endclaim

\proof. Let us consider the direct image $T$ of the current $\wh T$ on $X$.
It is a closed positive current in the class $c_1(L)$, and therefore the difference 
$$\wh T- \mu^\star T$$
lies in the class $\{\wh E\}$. On the other hand, we can write
$$\wh T= \chi_{\wh X\setminus E^\prime}\wh T+ \sum \rho^{(j)}E_j$$
and 
$$\mu^\star T= \chi_{\wh X\setminus E^\prime} \mu^\star T+ \sum \tau^{(j)}E_j$$
(where $E^\prime$ is the support of the exceptional divisor of the map $\mu$).
Since the two currents coincide outside $E^\prime$, we infer that
$$\wh T- \mu^\star T= \sum (\rho^{(j)}- \tau^{(j)})E_j.$$
We invoke now a cohomological argument: the classes associated to the hypersurfaces $E_j$ are independent in $H^{1,1}(X, \bR)$ and therefore
by the above relations we obtain $\rho^{(j)}- \tau^{(j)}= b^{(j)}$ for any $j$,
so our claim is proved.\hfill\qed

let us consider the space 
$${\cE}_{p_0}(u):= \{ U\in H^0(X, p_0(K_X+ S+ L)\big)
\hbox{ such that } U_{|S}= u\wedge (ds)^{m_0} \hbox { and } \star\star\}$$ 
where 
$$\int_{X} (U\wedge \ol U)^{1/m_0}\exp \bigl( -\varphi_L-\varphi_S\bigr)< \infty  \leqno (\star\star)$$
Remark in the first place that this space is not empty, by our previous considerations, and let us
take $U_{min}\in {\cE}_{m_0}(u)$ an element of our space with minimal $L^{2/m_0}$ norm.
Then we claim that $U_{min}$ is the extension we seek. 

Indeed, let us endow the bundle
$(m_0-1)(K_X+ S+ L)+ L$ with the metric given by the section $U_{min}$ raised to the power 
$1-1/m_0$, twisted with the metric $h_L$. By the Ohsawa-Takegoshi-Manivel theorem, there exist
$V\in H^0(X, m_0(K_X+ S+ L)\big)$ such that  $V_{|S}= u\wedge ds$ and such that 
$$\int_{X} {{V\wedge \ol V}\over {(U\wedge \ol U)^{{m_0-1}\over {m_0}} }}\exp \bigl( -\varphi_L-\varphi_S\bigr)\leq C_{univ}\int_S(u\wedge \ol u)^{1/m_0}\exp( -\varphi _L)$$
where $C_{univ}$ above is a numerical constant; remark that we use here the fact that $U_{min}$ is an extension of $u\wedge (ds)^{m_0}$. 

But now we will show that the minimality of $U_{min}$ imply 

$$\int_{X} (U_{min}\wedge \ol U_{min})^{1/m_0}\exp \bigl( -\varphi_L-\varphi_S\bigr)\leq 
\int_{X} {{V\wedge \ol V \exp \bigl( -\varphi_L-\varphi_S\bigr)}\over {(U_{min}\wedge \ol U_{min})^{{m_0-1}\over {m_0}} }}\leqno (15)
$$
and of course, once (15) is proved, we are done.
\vskip 10pt
In order to settle this last step we argue by contradiction: assume that (15) is false; then we have
$$\eqalign{
& \int_{X} (V\wedge \ol V)^{1/m_0}\exp \bigl( -\varphi_L-\varphi_S\bigr)= \cr 
& \int_X {{(V\wedge \ol V)^{{1}\over {m_0}} \exp \bigl( -{{\varphi_L+ \varphi_S}\over {m_0}}\bigr)}\over {(U_{min}\wedge \ol U_{min})^{{m_0-1}\over {m_0^2}} }}(U_{min}\wedge \ol U_{min})^{{m_0-1}\over {m_0^2}} \exp \big(
(1-1/m_0)(\varphi_L+ \varphi_S)
\big)\leq \cr
 & \Big(\int_{X} {{V\wedge \ol V \exp \bigl( -\varphi_L-\varphi_S\bigr)}\over {(U_{min}\wedge \ol U_{min})^{{m_0-1}\over {m_0}} }}\Big)^{1/m_0}
\Big(\int_{X} (U_{min}\wedge \ol U_{min})^{1/m_0}\exp \bigl( -\varphi_L-\varphi_S\bigr)\Big)^{1-1/m_0} \cr
 & < \int_{X} (U_{min}\wedge \ol U_{min})^{1/m_0}\exp \bigl( -\varphi_L-\varphi_S\bigr)\cr
}$$
Thus we see that our assumption (15) is valid 
and therefore the claim (?) is settled. \hfill \qed

\medskip

\noindent In this second part we  prove 
an extension theorem concerning  metrics with positive curvature current of adjoint 
$\bQ$--line bundles. Particular cases of the statement below (and its corollaries) were used several times in the 
proof of the main theorem. 

\noindent The general framework of our result is as follows. 

\noindent $\bullet$ We consider a projective, non-singular manifold $X$,
and let $\Delta_0\to X$ be a $\bQ$-line bundle, such that for any $m$ large and divisible enough, we have
$$m\Delta_0\equiv L_1+\cdots L_{m-1}\leqno (37)$$
where for each $1\leq j\leq m-1$ the $L_j$ is a line bundle, endowed with a 
metric $h_j$ whose curvature current is assumed to be positive. 
\smallskip
\noindent $\bullet$ $\Delta_1\to X$ is a $\bQ$--line bundle, whose first Chern class
contains a closed positive current $\Xi$. We introduce the  notation
$$\Delta:= \Delta_0+ \Delta_1.$$
\smallskip

\noindent $\bullet$  $S\subset X$ is a non-singular hypersurface such that $h_{j|S}$ 
and $\Xi_{|S}$ are well-defined,
in the sense that the local potentials of $h_j$ are not identically $-\infty$ when restricted to $S$. We assume the existence of a closed positive current 
$$T\in c_1(K_X+ S+ \Delta)$$
such that:
{\itemindent 6mm

\smallskip
\item {$(\cR)$} The restriction of $T$ to $S$ is well defined; 
\smallskip
\item {$(\cT)$} There is some positive number $\varepsilon$ such that for any $\delta>0$ and  any $j= 1,..., m-1$ we have 
$$\cI \big( (1-\delta)\varphi_{L_j}+ \varepsilon \varphi_{T|S}\big)= \cO_S.\leqno (\cT_j)$$
}
\medskip

\noindent We observe that under the hypothesis $(\cT_j)$ the multiplier sheaf
of the metric $(1-\delta)\varphi_j$ is automatically trivial for any $\delta> 0$ . The assumption above
intutively means  that the singularities of $\varphi_j$ and $\varphi_T$ are transversal.
\bigskip
The following statement is a {\sl metric version} of the extension theorem due to
Ein-Popa and Hacon-McKernan [14], [16] (as it will be seen in the corollaries which  follow  the proof). The origin of  this kind of results is the seminal article of Y.-T. Siu  [32].

\claim Theorem 2.1|Let $S$ be a non-singular hypersurface of a projective manifold $X$
and let $\Delta_0, \Delta_1, T$ be respectively $\bQ$--line bundles and a closed positive current with the properties stated above.
We consider a metric $h_{K_S+ \Delta}$ on the $\bQ$-line bundle
$K_S+ \Delta$ such that:
{\itemindent 5mm

\smallskip
\item {i)} The corresponding curvature current is positive;
\smallskip
\item {ii)} We have
$$\varphi_{K_S+ \Delta}\leq \varphi_{\Xi|S},$$ 
or more generally, it is enough to assume that for any positive integer $m$ we have
$$\cI (m\varphi_{K_S+ \Delta})\subset \cI(m\varphi_{\Xi|S}).\leqno (26)$$
}
\noindent Then there exist an ample line bundle $A\to X$ such that 
for any positive 
integer $m$,  
the metric $\displaystyle h_{K_S+ \Delta}$ admits a sub-extension 
$\displaystyle h^{(m)}$, where  $\displaystyle h^{(m)}$ is a metric of the $\bQ$--line bundle
$K_X+ S+ \Delta+1/mA$ defined on $X$, with positive curvature current.\hfill \qed
\endclaim
\medskip
\noindent In other words, we claim the existence of a metric  
$\displaystyle h^{(m)}$ of the $\bQ$--line bundle
$K_X+ S+ \Delta+1/mA$ with positive curvature and such that its restriction to $S$ is well defined and satisfies the next inequality
$$\varphi_{K_S+ \Delta}\leq \varphi^{(m)}- \varphi_A/m,$$ 
at each point of $S$, where $\varphi_A$ is a smooth metric on $A$.

\proof. Our arguments involve two main ingredients: an approximation theorem for  closed positive currents (which is originally due to Demailly in [9], and further polished by Boucksom in [6]) and the invariance of plurigenera  techniques. The first technical tool allows 
us to ``convert" the metric $\displaystyle h_{K_S+ \Delta}$ to a family of sections of the bundle $m(K_S+ \Delta)+ A$, whereas the second will show that each member of the
family of sections obtained in this way admits an extension to $X$. We remark that for the second step we could just quote
the recent result of Ein-Popa [14]; however, we prefer to present an alternative proof, which is more analytic in spirit (the reader can profitably consult the original proof in [14], based on the 
theory of {\sl adjoint ideals}).
\medskip

\subsection {\S 2.1 Approximation of closed positive currents}
\smallskip
\noindent To start with, we  recall the (crucial) fact that a closed positive current of 
$(1, 1)$--type
can be approximate in a very accurate way by 
a sequence of currents given by algebraic metrics.

\claim Theorem ([6], [9])|There exist an ample line bundle $A \to X$
and a smooth metric $h_A$ on $A$ 
such that the following holds. Let $m$ be a large and divisible enough positive integer, and let $\big(\sigma_j^{(m)}\big)$ be an orthonormal basis of the space
$$H^0\big(S, m(K_S+ \Delta)+ A\big)$$
endowed with the metric whose local weights are given by the functions
$\displaystyle m\varphi_{K_S+ \Delta}+ \varphi_{A}$ and an arbitrary smooth volume form on $S$. Then we have the pointwise inequality
$$\varphi_A/m+\varphi_{K_S+ \Delta}\leq {1\over m}\log\big(\sum_{j=1}^{N_m}|f_j^{(m)}|^2\big)+ C\leqno (27)
$$
where $f_j^{(m)}$ is the local holomorphic function corresponding to the 
section $\sigma_j^{(m)}$, for any $j$.\hfill \qed

\endclaim

\medskip
\noindent We remark that the precise statement in [6], [9] is slightly different, 
but the proof given there clearly contains the inequality (27) above. 
Notice that if the metric $\varphi_{K_S+ \Delta}$ satisfies condition (ii) of Theorem 2.1, then 
the sections $(\sigma_j^{(m)})$ satisfy the integrability relation 
$$\sigma_j^{(m)}\in \cI (m\varphi_{\Xi |S})\leqno (28)$$
for any positive integers $j, m$.
\medskip
\subsection {\S 2.2 Extension of sections}
\smallskip

\noindent Let $u$ be a section of $m(K_S+\Delta)+A$
such that 
$$\int_S\vert u\vert ^2\exp (-m\varphi_{K_S+ \Delta}-\varphi_A)d\lambda< \infty;\leqno (29)$$
so that in particular $u\in \cI (m\varphi_{\Xi |S})$ (for example, $u$ may be 
one of the $\sigma_j^{(m)}$ involved in the approximation process above).

The main part of the proof of 2.1 is to show that $u$ extends to $X$.
(Of course, we are going to use heavily the other metric assumptions 
in the statement 2.1, i.e. just the relation (29) above is  not sufficient). 
To this end, we will use the invariance of plurigenera techniques. 
The argument is derived from the original article
of Siu, see [32], as well as [13], [26].

\noindent We denote by $L_m:= m\Delta_1$ and
$$L^{(p)}:= p(K_X+ S)+L_1+...+ L_p\leqno (30)$$
where $p= 1,..., m$. By convention, $L^{(0)}$ is the trivial bundle.
We can assume that the bundle $A\to X$ in 
the approximation statement above is ample 
enough, but independent of $m$, so that the next two conditions are satisfied.

{
\itemindent 6mm
\noindent
\item {$(\dagger)$} For each $0 \leq p\leq m-1$, the bundle
$L^{(p)}+ mA$ is generated by its global sections, which we denote by
$(s^{(p)}_j)$;
\smallskip
\item {$(\dagger\dagger)$} Any section of the bundle $L^{(m)}+
mA_{\vert S}$ admits an extension to 
$X$.

}
\medskip
\noindent  
We formulate now the following inductive statement (depending on a positive integer $k$ and 
on $0\leq p\leq m-1$):
\vskip 7pt
\noindent 
{$\bigl({\cal P}_{k,p}\bigr):$ 
\sl The sections 
$$u^k\otimes s^{(p)}_j\in H^0\bigl(S, L^{(p)}+ kL^{(m)}+ (k+m)A_{\vert S}\bigr)$$
extend to $X$, for each $j= 1,..., N_p.$}
\vskip 7pt
\noindent If we can show that ${\cal P}_{k,p}$ is true for any $k$ and $p$, 
then we are done, since the extensions $U^{(km)}_j$ of the sections 
$u^k\otimes s^{(0)}_j$ can be used to define a metric
$$
\log\sum|U_j^{(km)}|^2
$$
on the bundle
$$km(K_X+S+ \Delta)+ (k+m)A$$
whose $km^{\rm th}$ root it is defined to be 
$h^{(m)}$ (for $k\gg 0$). Theorem A shows that  $h^{(m)}$
is then indeed a sub-extension of the given metric $h_{K_S+ \Delta}$.

\medskip
\noindent Thus, it is enough to check the property ${\cal P}_{k,p}$; our
arguments rely heavily on the next version of the
Ohsawa-Takegoshi theorem, proved by L. Manivel (see [1], [2], [12], [22], [24], [25] to quote
only a few), and further refined by McNeal-Varolin in [23].

\claim 2.2 Theorem ([23])|Let $X$ be a projective
$n$-dimensional 
manifold, and let $S\subset X$ be a non-singular hypersurface.
Let $F$ be a line bundle,
endowed with a metric $h_F$. We assume the existence 
of some non-singular metric $h_S$ on the bundle $\cal O(S)$ such that:  

{\itemindent 6mm
\smallskip
\item {(1)} $\displaystyle {{\sqrt {-1}}\over {2\pi}}\Theta_F
\geq 0$ on $X$;
\smallskip
\item {(2)} $\displaystyle {{\sqrt {-1}}\over {2\pi}}\Theta_F -  
\alpha
 {{\sqrt {-1}}\over {2\pi}}\Theta_S\geq 0$ 
for some $\alpha> 0$;
\smallskip
\item {(3)} The restriction
of the metric $h_F$ on $S$ is well defined.

}
\smallskip

\noindent Then every section $u\in H^0\bigl(S, (K_X + S+
F_{\vert S})\otimes {\cal I}(h_{F|S})\bigr)$ admits an extension
$U$ to $X$ such that 
$$c_n\int_XU\wedge \ol U\exp \Big(-\varphi_F-\varphi_S-\log\big(|s|^2\log^2(\vert s\vert)\big)\Big)< \infty$$
where $s$ is a section whose zero set is precisely the hypersurface $S$ and its norm in the integral above is measured with respect to $h_S$. In particular
$$
c_n\int_XU\wedge \ol U\exp \Big(-\varphi_F-\varphi_S\Big)/|s|^{2(1-\delta)}< \infty$$
for any $\delta>0$.
\endclaim
\medskip 
\claim 2.3 Remark|{\rm In fact, the general result in [23] allow to extend sections 
of bundles $K_X+ E_{|W}$ defined on higher codimensional manifolds $W\subset X$
which are complete intersections 
$$W= \cap_j (\rho_j= 0)$$
where $\rho_j$ are sections of some line bundle $\Lambda$. Of course, in this case the curvature and integrability conditions $(1)-(3)$ are modified: basically we replace 
$\displaystyle {{\sqrt {-1}}\over {2\pi}}\Theta_F$ with 
$\displaystyle {{\sqrt {-1}}\over {2\pi}}\Theta_E+ \displaystyle {{\sqrt {-1}}\over {2\pi}}
\ddbar \log \big(\sum |\rho_j|^2\big)$, and the curvature term of $\cO(S)$ with 
the curvature of $\Lambda$  
(if $W= S$ is a hypersurface, the bundle $E$ in this setting 
is just $F+ S$). We refer to [23] for the precise statement.\hfill \qed
}
\endclaim

\medskip
We will use inductively the theorem above, in order to show the extension property
${\cal P}_{k,p}$. The first steps are as follows.

\noindent (1)  For each $j= 1,...,N_0$, the section 
$u\otimes s^{(0)}_j\in H^0\bigl(S, L^{(m)}+ (m+1)A_{\vert S}\bigr)$
admits an extension $U^{(m)}_j\in H^0\bigl(X,  L^{(m)}+ (m+1)A\bigr)$,
by the property $\dagger\dagger$.

\noindent (2) We use the sections $(U^{(m)}_j)$ to construct a metric 
$$\varphi^{(m)}=\log\sum|U_j^{(m)}|^2$$
on the bundle $L^{(m)}+ (m+1)A$.

\noindent (3) Let us consider the section 
$u\otimes s^{(1)}_j\in H^0\bigl(S, L^{(1)}+ L^{(m)}+ (1+m)A_{\vert S}\bigr)$. 
We remark that the bundle
$$L^{(1)}+ L^{(m)}+ (m+1)A = K_X+ S+ L_1+ L^{(m)}+ (m+1)A$$
can be written as $K_X+ S+ F$ where
$$F:=  L_1+ L^{(m)}+ (m+1)A.$$
We are going to construct a metric on $F$ which satisfy the 
curvature and integrability assumptions in the theorem 2.2.

Let $0< \delta\ll \varepsilon\ll 1$ be positive real numbers. We endow the bundle 
$F$ with the metric given by
$$\varphi^{(m)}_{\delta, \varepsilon}:= (1-\delta)\varphi_{L_1}+ \delta\wt \varphi_{L_1}+ 
(1-\varepsilon)\varphi^{(m)}+ \varepsilon((m+1)\varphi_A+ m\varphi_T)\leqno (31)$$
where the metric $\wt \varphi_{L_1}$ is smooth (no curvature requirements) and  
$\varphi_{L_1}$ is the singular metric whose existence is given by hypothesis.
The curvature conditions (1) and (2) in 2.2 are satisfied, since we can use the metric on $A$
in order to absorb the negativity of $\wt \varphi_{L_1}$. 

Next we claim that the sections $u\otimes s^{(1)}_j$ are integrable with respect to the metric defined in (31).
Indeed, the integrand is less singular than 
$$
\exp \big(-(1-\delta)\varphi_{L_1}- \varepsilon\varphi_T\big)$$
which is integrable  by the hypothesis $(\cT_1)$.
\smallskip

\noindent (4) We apply the extension theorem 2.2 and we get
$U^{(m+ 1)}_j$, whose
restriction on $S$ is precisely
$u\otimes s^{(1)}_j$.\hfill \qed
\medskip

\vskip 5pt Now the assertion  ${\cal P}_{k,p}$ will be obtained by iterating
the
procedure (1)-(4) several times. Indeed, assume that the 
proposition ${\cal P}_{k,p}$ has been checked, and consider the set of global sections 
$$U^{(km+p)}_j\in H^0\bigl(X, L^{(p)}+ kL^{(m)}+ (k+m)A\bigr)$$
which extend $u^k\otimes s^{(p)}_j$. They induce a metric on the above bundle, denoted by $\varphi^{(km+p)}$. 

\noindent If $p< m-1$, then we define the family of 
sections 
$$u^k\otimes s^{(p+1)}_j\in H^0(S, L^{(p+1)}+ kL^{(m)}+ (k+m)A_{|S})$$
on $S$. 
As in the step (3) above  we have that
$$L^{(p+1)}= K_X+ S+ L_{p+1}+ L^{(p)}.$$
To apply the extension result 2.2, we have to exhibit a metric
on the bundle 
$$F:=  L_{p+1}+ L^{(p)}+kL^{(m)}+ (k+ m)A$$
for which the curvature conditions are satisfied, and such that the family of sections above are $L^2$ with respect to it.
We define
$$\varphi^{(km+p+1)}_{\delta, \varepsilon}:= (1-\delta)\varphi_{L_{p+1}}+ \delta\wt \varphi_{L_{p+1}}+ 
(1-\varepsilon)\varphi^{(km+ p)}+ \varepsilon( \wt \varphi_{L^{(p)}}+ mk\varphi_T+ (m+k)\varphi_A)\leqno (32)$$
and we see that the parameters $\varepsilon, \delta$ have to satisfy the conditions:

{\itemindent 6mm

\item {(A)} We don't assume any relation between the zero set of $u$ and the singularities of $T$, thus we have to keep the poles of $T$ ``small" in the expression of the metric above, as we have to use 
$(\cT_{p+1}$). This impose $mk\varepsilon\ll 1$;

\item {(B)} We have to absorb the negativity in the smooth curvature terms in (32), thus we are forced to have $\delta\ll\varepsilon$.

}
\noindent Clearly we can choose $\varepsilon, \delta$ with the properties required above. Let us check next the $L^2$ condition. As before, the integrand is less singular than
$$
\exp \big(-(1-\delta)\varphi_{L_{p+1}}- mk\varepsilon\varphi_T\big).$$
This is again integrable  because of the transversality hypothesis $(\cT_{p+1}$).

\medskip
Finally, let us indicate how to perform the induction step if $p= m-1$:
we consider the family of 
sections 
$$u^{k+1}\otimes s^{(0)}_j\in H^0(S, (k+1)L^{(m)}+ (m+k+1)A_{|S}),$$ 
In this case we have to exhibit a metric
on the bundle 
$$L_{m}+ L^{(m-1)}+kL^{(m)}+  (m+k+1)A.$$
This is easier than before, since we can simply take 
$$\varphi^{m(k+1)}:= 
 m\varphi_{\Xi}+ \varphi_A+ 
\varphi^{(km+m-1)}\leqno (33)$$
With this choice, the curvature conditions are satisfied; as for the $L^2$ ones, we remark that we have
$$\eqalign{
 & \int_S{{|u^{k+1}\otimes s^{(0)}_j|^2}\over {(\sum_q \vert u^k\otimes s^{(m-1)}_q\vert^2)}}
\exp \big(-m\varphi_{\Xi}\big)dV< \cr
 < & C\int_S{{|u\otimes s^{(0)}_j|^2}\over {(\sum_q \vert s^{(m-1)}_q\vert^2)}}
\exp \big(-m\varphi_{\Xi}\big)dV< \infty \cr
}$$
where the last relation holds because of (29).
The proof of the extension theorem is therefore finished.\hfill \qed
\vskip 10pt

\bigskip
\noindent Next, we will highlight 
a very useful geometric context for which the 
numerous hypothesis in the theorem 2.1 are satisfied.

\noindent To this end, let us consider the next objects:
{\itemindent 5mm
\smallskip
\item {$(1)$} $\Delta_0$ is an effective $\bQ$-divisor and has critical exponent greater than 1
(i.e. it is klt);
\smallskip
\item {(2)} $\Delta_1, \Delta _2$ are effective, respectively ample $\bQ$--bundles;
\smallskip
\item {(3)} The support of $S+ \Delta_0+ \Delta_1$ has normal crossing
and $S$ is not contained in the support of $\Delta_0+ \Delta_1$;
\smallskip
\item {(4)} Let $\Delta:= \Delta_0+ \Delta_1+\Delta_2$; we assume that there exist a closed positive current 
$$T\in c_1(K_X+ S+ \Delta)$$
such that the restriction of $T$ to any intersection of components of the support of the divisor $S+ \Delta_0+ \Delta_1$ is well defined (compare with the Hacon-McKernan
assumption on the stable base loci in [16]).

}

\medskip 
\noindent Our claim is that the assumptions in the theorem 2.1 are satisfied, 
provided that $(1)-(4)$ are. Indeed, let us first define the line bundles $L_j$
such that 
$$m\Delta_0= L_1+\cdots + L_{m-1}$$
and the corresponding metrics on them.

 By hypothesis, the $\bQ$--divisor $\Delta_0$ is effective, its critical exponent is
greater than one and the divisors in his support have normal crossing. Therefore, we have the decomposition
$$\Delta_0= \sum_{j= 1}^N{{p_j}\over {q_j}}Z_j$$
where the $(p_j, q_j)$ above are positive integers, such that coefficients above are {\sl strictly smaller   
than 1}; we will assume that the sequence $\big({{p_j}\over {q_j}} \big)$ is increasing. 
As a consequence of the above decomposition, we see that for any positive and divisible enough 
$m$ we can find the line bundles 
$L_1,..., L_{m- 1}$ such that 
$$m\Delta_0\equiv L_1+...+ L_{m-1}$$
and such that 
$$L_j\equiv \sum_{p\in \Lambda_j} Z_p$$
with $\Lambda_j\subset \{1,...,N\}$. To see this, we first remark that 
for each index $j$ we have $\displaystyle m_j:= m{{p_j}\over {q_j}}\leq m-1$, and then 
we define 
$$\Lambda_j:= \{1,..., N\}$$ 
for $j= 1,...,m_1$. The next package of $(L_j)$ is defined 
by the set 
$$\Lambda_j:= \{2,...N\}$$
for $j:= m_1+1,..., m_2$, and so on. In this way we have our bundles 
$L_1,...L_{m-1}$; we endow the corresponding $L_j$
with the natural singular metric induced by the hypersurfaces $(Z_p)$
whose indexes belong to $\Lambda_j$.  
We define 
$$L_m:= m(\Delta_1+ \Delta_2).$$
In view of the assumptions $(1)-(4)$ above,
all that still needs to be checked is the fact that the multiplier sheaves $(\cT_j)$ are trivial; we are going to explain this along the following lines.

\bigskip

We remark next that the setting $(1)-(4)$ is precisely the context in which we have
found ourselves in the section 1.3. We recall now the result by Ein-Popa (see [14], as well as [7], [13], [26], [27], [29], [31], [32], [33], [34]).

\claim 2.3 Theorem {\rm ([14])}|Let $X$ be a non-singular projective manifold. Let $S\subset X$ be a non-singular hypersurface,
and let $\Delta= \Delta_0+\Delta_1+\Delta_2$ be a sum of $\bQ$--line bundles on $X$
which verify the properties 1)-4) above.
Then any section $u$ of the bundle 
 $m(K_{X}+ S+ \Delta)_{|S}$ whose zero set contains the divisor $m\Delta_1$
 extends to $X$, for any $m$ divisible enough. \hfill \qed
 \endclaim
 
\proof. We will use the notations in the proof of the theorem 2.1; it was established
during this proof that the family of sections $\big(u^k\otimes s^{(0)}_j\big)$
 admits an extension to $X$, which was denoted by
 $$U^{(km)}_j\in H^0\bigl(X, kL^{(m)}+ mA\bigr)$$
 Here we remark that we have the factor $mA$ instead of $(k+m)A$
 in the 2.1 above, simply because $u$ is a section of the bundle 
 $m(K_{X}+ S+ \Delta)_{|S}$ (and not of $A+ m(K_{X}+ S+ \Delta)_{|S}$
 as in 2.1).
We take $k\gg 0$, so that
$$ {{m-1}\over {k}}A < \Delta_2,
$$
in the sense that the difference is
ample (we recall that $\Delta_2$ is ample). Let $h_A$ be a smooth, 
positively curved metric on $\Delta_2-{{m-1}\over{k}}A$. 
\vskip7pt
\noindent We apply now the extension theorem 2.2 with the following data
$$\displaystyle F:= {{m-1}\over {m}}L^{(m)}+ \Delta_0+ \Delta_1+ \Delta_2$$
and the metric $h$ on $F$ constructed as follows:
\smallskip
\item {(1)} on the factor  $\displaystyle {{m-1}\over {m}}L^{(m)}+ {{m-1}\over{k}}A$
we consider the appropriate power of the metric $\varphi^{(km)}$
given by the sections $(U^{(k)}_j)$; on $\Delta_0$ and $\Delta_1$ we consider the 
singular metric induced by the corresponding divisors, and finally on 
$ \Delta_2- {{m-1}\over{k}}A$ we take the metric $h_A$;
 
\smallskip
\item{(2)} we take an arbitrary, smooth metric $h_S$
on the bundle associated to the hypersurface $S$.
\smallskip
\noindent With our choice of the bundle $F$, 
the section $u$ we want to extend become a section of the adjoint 
bundle $K_X+ S+F_{\vert S}$ and the positivity requirements (1) and (2) in the extension
theorem are satisfied since we have
$$ {{\sqrt {-1}}\over {2\pi}}\Theta_h(F) \geq {{\sqrt {-1}}\over {2\pi}}\Theta_{h_A}(\Delta_2- {{m-1}\over
{mk}}B)> 0$$

Concerning the integrability of $u$,
remark that we have
$$\eqalign{
\int_S{{\vert u\vert^2\exp(- \varphi_{\Delta_10}- \varphi_{\Delta_1}-
\varphi_{A})}\over 
{\bigl(\sum_j\vert U^{(km)}\vert^2\bigr)^{{m-1}\over {mk}}}}= & 
\int_S{{\vert u\vert^2\exp(- \varphi_{\Delta_0}- \varphi_{\Delta_1}-
\varphi_{A})}\over 
{\bigl(\sum_j\vert u^{\otimes km}\otimes
s_j^{(0)}\vert^2\bigr)^{{m-1}\over {mk}}}}\leq \cr
\leq &C \int_S\vert u\vert^{2/m}\exp(-\varphi_{\Delta_1}-
\varphi_{A})\exp (- \varphi_{\Delta_0})dV<\cr
< & \infty\cr
}
$$
and this last integral converge by the hypothesis concerning the 
zero set of $u$ and the fact that the 
multiplier ideal sheaf of the restriction of $\Delta_0$ to $S$
is trivial.\hfill \qed
\medskip


We recall the notations that have been used in the introduction.
Let $L\to \cX$ be a line bundle, endowed with a non-singular hermitian metric
$h_{L, 0}$. We assume the existence of a function $\psi\in L^1(\cX)$ such that 
the curvature current of $L$ with respect to the metric 
$h_L:= \exp(-\psi)h_{L, 0}$ is positive, moreover, let us suppose that the restriction of $\psi$
to the central fiber $\cX_0$ is not identically 
$-\infty$. 

Over the central fiber $\cX_0$, we have closed positive current
$$T_0:= m_0\Theta_h(K_{{\cal X}_0})+ \Theta_{h_{L,0}}(L)+ \sqrt {-1}\ddbar \varphi_0$$
and by hypothesis we assume that $\sup_{\cX_0}(\varphi_0- \psi)< \infty.$

Let $(A, h_A)$ be a line bundle endowed with a non-singular metric whose curvature 
form (denoted by $\omega$) is positive. For each positive integer $m$, we will consider the space 
$$\cV_m\subset 
H^0\bigr( \cX_0, m(m_0K_{\cX_0}+ L)+A\bigl) $$ 
which consists of sections  $u$ with finite $L^2$ norm
$$||u||^2:= \int_{\cX_0}\vert u\vert ^2\exp (-m\varphi_{T_0}-\varphi_A)dV_\omega$$
(the notation above just means that the section $u$ is $L^2$ integrable with respect to
the metric induced by the current $T_0$). 

We consider a orthonormal basis of the space above, denoted by $(u^{(m)}_j)$,
and let $T_0^{(m)}$ be the curvature current of the metric associated to the sections 
above, normalized by the factor $1/m$. Thus the precise expression of this current is
$$T_0^{(m)}= m_0\Theta_h(K_{{\cal X}_0})+ \Theta_{h_{L,0}}(L)+ {1\over m}\Theta_{h_A}(A)
+ \sqrt {-1}/m\ddbar \log \bigl( \sum_j|u^{(m)}_j|^2_{h^{mm_0}h_{L,0}^mh_A}\bigr).$$

\noindent Then we recall the next statement.

\claim Theorem {\rm [7], [12]}|The current $T_0^{(m)}$ converge weakly to $T_0$, as 
$m\to\infty$. Moreover, its potential converge pointwise (and in $L^1$ norm) to
the one of $T_0$.
\endclaim

\noindent In other words, the function

$$\varphi^{(m)}_0:= 1/m\log \bigl( \sum_j|u^{(m)}_j|^2_{h^{mm_0}h_{L,0}^mh_A}\bigr)$$converge to $\varphi_0$ at each point of 
$\cX_0$. 

We come now to the main part of the proof, namely the extension of each
section $u^{(m)}_j$ with precise estimates.

\bigskip
\noindent In connection with the theorem we have just proved, we would like to formulate the
following problem.

\claim 4.2.3 Conjecture|Let $X$ be a non-singular projective manifold and let 
$S\subset X$ be a smooth hypersurface. Consider 
$L\to X$ be a line bundle, endowed with a metric $h_L$ such that 
the curvature current verify the condition 
$$\sqrt {-1} \Theta_{h_L}(L)\geq \varepsilon_0\sqrt {-1}\Theta_{h_A} (A)$$
for some positive $\varepsilon$ and some non-singular metric with 
positive curvature on the ample bundle $A$. Let
$$T_S:= \sum_{j=1}^N\nu_j[\Xi_j]$$
be a closed positive current on $S$, where $\nu_j\in \bR_+$ and the $\Xi_j\subset S$ are irreducible 
and reduced hypersurfaces; moreover we assume that 
$$T_S\in c_1(K_X+ S+ L_{|S})$$
and that the local potentials of the current $T_S$ are bounded with respect to the restriction of the 
metric of $L$ to $S$.

Then there exist a closed positive current
$$T\in c_1(K_X+ S+ L)$$
with the following properties
{\itemindent 4mm
\smallskip
\item {i)} The restriction of $T$ to $S$ is well defined, and it is less singular than $T_S$;
\smallskip
\item {ii)} $T$ is the current of integration over an effective $\bR$-divisor on $X$.

}
\endclaim

\noindent Instead of $ii)$ above, 
a more realistic expectation would be that $T$ has analytic singularities i.e. locally its potentials are 
equal to the logarithm of the sum of squares of holomorphic functions, modulo smooth functions.

This question is motivated by the work of Y.-T. Siu  [37] concerning the finite generation property of the canonical ring of the algebraic varieties of general type. 
\vskip 7pt
\claim 4.2.4 Remark|{\rm Under the hypothesis of the above conjecture, the proof of the 0.2 
show the existence of 
a closed positive current
$$T\in c_1(K_X+ S+ L)$$
whose restriction to $S$ is less singular than $T_S$; this is due to the fact that the elements of the 
space $\cV_m$ above are {\sl grosso modo} given by powers of the sections associated to the hypersurfaces $\Xi_j$. Of course, since we are in the "compact" setting, we have to use the appropriate version of the extension theorem of Ohsawa-Takegoshi type in the above proof, 
see [15].
Unfortunately, even in the case where the line bundle $L$ is big, we have to take the limit
in order to remove the auxiliary line bundle $A$ and this process destroy the analytic singularities 
of the currents we get during the inductive process. 
}
\endclaim

On the other hand, according to a recent result of Campana-Peternell (see [9], theorem 3.1)
{\sl if the canonical class of a projective manifold $X$ contains an effective $\bQ$--divisor, then some multiple of $K_X$ has a section}. This result is a rather strong support for the above conjecture.

The results in the field dealing with this situation (see [22], [16], [17], [18]
to quote only a few) are assume that

\noindent Many interesting results having recently emerged from 
algebraic
geometry ([20], [44], [5] to quote only a few) are dealing with 
extension of sections of bundles like 
$$pK_{\cX}+ L_{|\cX_0}\leqno(3)$$ where $\cX_0$ is the central fiber of a projective family $\pi: \cX\to \bD$ over the unit disk,
and $(L, h_L)$ is a positively curved 
line bundle, such that the multiplier ideal sheaf of the restriction metric
$\displaystyle h_{L|\cX_0}^{1/p}$ is equal to the structural sheaf. 

Given the results in the literature, we remark the following fact.
If the Chern class of $L$ is divisible by $p$, then any section $u$ of the bundle above extends, by the results of B. Claudon, J.-P. Demailly, S. Takayama, D. Varolin in . However, if this is not the case, then the articles of Hacon-McKernan, Ein-Popa, Hacon-de Fernex () prove that 
additional {\sl vanishing hypothesis} for the section together with {\sl strict positivity} of the curvature current associated to $(L, h_L)$ are {\sl sufficient}
in order to show that $u$ extends. In view of the examples in [22], we see that {\sl not all} of these assumptions can be removed.